\documentclass[a4paper,10pt]{article}
\usepackage{stmaryrd}
\usepackage{amsfonts}
\usepackage{bbm}
\usepackage{amscd}
\usepackage{mathrsfs}
\usepackage{latexsym,amssymb,amsmath,amscd,amscd,amsthm,amsxtra,xypic}
\usepackage[dvips]{graphicx}
\usepackage[all]{xy}
\usepackage{ulem}

\newtheorem{thm}{Theorem}[section]
\newtheorem{lem}[thm]{Lemma}
\newtheorem{cor}[thm]{Corollary}
\newtheorem{pro}[thm]{Proposition}
\newtheorem{ex}[thm]{Example}
\newtheorem{rmk}[thm]{Remark}
\newtheorem{defi}[thm]{Definition}

\newcommand{\be }{\begin{eqnarray*}}
\newcommand{\ee }{\end{eqnarray*}}

\newcommand{\pf}{\noindent{\bf Proof.}\ }

\newcommand{\frkg}{\mathfrak g}

\newcommand{\g}{\mathfrak g}
\newcommand{\m}{\mathfrak m}

\def\gpd{\,\lower1pt\hbox{$\longrightarrow$}\hskip-.24in\raise2pt
         \hbox{$\longrightarrow$}\,}

\def\qed{\hfill ~\vrule height6pt width6pt depth0pt}

\newcommand{\Courant}[1]{\left\llbracket  #1\right\rrbracket }

\newcommand{\Hom}{\mathrm{Hom}}
\newcommand{\Der}{\mathrm{Der}}

\newcommand{\Aut}{\mathrm{Aut}}
\newcommand{\gl}{\mathrm {gl}}

\newcommand{\Img}{\mathrm{Im}}

\newcommand{\End}{\mathrm{End}}
\newcommand{\Inn}{\mathrm{Inn}}
\newcommand{\DER}{\mathrm{DER}}
\newcommand{\Ha}{\mathrm{H}}
\newcommand{\coker}{\mathrm{coker}}

\newcommand{\degree}{\mathrm{degree}}
\newcommand{\La}{\mathrm{L}}
\newcommand{\D}{\mathrm{D}}

\begin{document}
\title{
Crossed modules for Lie $2$-algebras}
\author{{ Honglei Lang and Zhangju Liu} \\
{Department of Mathematics and LMAM}\\ {Peking University,
Beijing 100871, China}\\
{\sf email: hllang@pku.edu.cn; ~~ liuzj@pku.edu.cn} }


\date{}

\maketitle
\begin{abstract}
The notion of crossed modules for Lie $2$-algebras is introduced. We
show that, associated to such a crossed module,  there is a strict
Lie $3$-algebra structure on its mapping cone complex and a strict
Lie $2$-algebra structure on its derivations. Finally, we classify
strong crossed modules by means of the third cohomology group of Lie
$2$-algebras.

\end{abstract}
\tableofcontents
\section{Introduction}
Crossed modules of Lie algebras first appeared in the work of
Gerstenhaber (\cite{Gerstenhaber}), which can be classified by  use of
the third cohomology group  of Lie algebras as follows: for a
crossed module of Lie algebras $\varphi:\m\longrightarrow{\g}$,
there exists a four term exact sequence
\begin{equation*}
\xymatrix@C=0.5cm{ 0 \ar[r] & \mathbb{V }\ar[rr]^{i} &&
                \mathfrak{m} \ar[rr]^{\varphi} && \mathfrak{g} \ar[rr]^{\pi} && \mathfrak{h}\ar[r]  & 0,
                }
\end{equation*}
where the cokernel $\mathfrak{h}$ is a Lie algebra and the kernel
$\mathbb{V}$ is an $\mathfrak{h}$-module induced by the action of
$\g$ on $\m$. Denote by crmod$(\mathfrak{h},\mathbb{V})$ the set of
equivalence classes of crossed modules with fixed kernel
$\mathbb{V}$, cokernel $\mathfrak{h}$ and action. Gerstenhaber
proved that there is a bijection between
crmod$(\mathfrak{h},\mathbb{V})$ and
$\Ha^3(\mathfrak{h},\mathbb{V}).$ See also \cite{Wagemann} for more
details and \cite{Baues,Casas} for other algebraic structures.

Lie algebras can be categorified to Lie 2-algebras.  For a good
introduction on this subject see \cite{Baez,Lada,Zhuo}. A  Lie $2$-algebra is a $2$-vector space equipped with
 a skew-symmetric bilinear functor, such that the Jacobi identity is controlled by a natural isomorphism,
which satisfies the coherence law of its own.
It is well-known that the notion of  Lie $2$-algebras is  equivalent
to that of 2-term $L_\infty$-algebras and the category of strict Lie 2-algebras is isomorphic to the category of crossed modules of Lie algebras (\cite{Baez}).

The cohomology of $L_\infty$-algebras and $A_\infty$-algebras was
studied in \cite{Lada,Penkava,Markl1} and in \cite{AQ} for a more
general theory.  While the cohomology of Lie 2-algebras was
formulated in \cite{Sheng2} for the strict case to characterize
strict Lie 2-bialgebras and in \cite{Liu} for the general case to
depict the deformation of Lie 2-algebras. For the cohomology of Lie
$2$-groups, see \cite{Ping}. In this paper, we propose the notion of
crossed modules of Lie $2$-algebras (Definition \ref{defi:Lie 2 cm})
and classify strong crossed modules via the third cohomology group
of Lie $2$-algebras (Theorem \ref{thm:main}).



Moreover,
given a crossed module of Lie $2$-algebras
$(\m,\g,\phi,\varphi,\sigma)$, there is a strict Lie $3$-algebra
($l_4=0$) structure on its mapping cone complex:
$\m_1\longrightarrow\g_1\oplus \m_0\longrightarrow \g_0$ (Theorem
\ref{thm:Lie 3}). Also, we obtain a strict Lie $2$-algebra structure
on derivations
$\Der(\g,\m):\Hom(\g_0,\m_1)\stackrel{-\D}{\longrightarrow}\Der_0(\g,\m)$
(Theorem \ref{thm:Der(g,m)}), where $\Der_0(\g,\m)$ is the set of
$1$-cocycles  and $\D$ is the Lie $2$-algebra coboundary operator.
Moreover, we get a Lie algebra structure on the first cohomology
group $\Ha^1(\g,\m)$.

This paper is organized as follows: In Section \ref{section2}, we
sketch some background on Lie 2-algebras, including basic
definitions, the cohomology theory and the derivations of Lie
$2$-algebras. In Section \ref{section cm}, we introduce the notion
of crossed modules of Lie $2$-algebras with some examples and
demonstrate that there is a strict Lie 3-algebra on the mapping cone
complex. In Section \ref{1-degree}, we provide a Lie algebra
structure on the set of $1$-cocycles. Then we prove that for a
crossed module, there exists a strict Lie $2$-algebra structure on
its derivations $\Der(\g,\m)$ and a Lie algebra structure on
$\Ha^1(\g,\m)$. Section \ref{section6} is concerned about the
classification of strong crossed modules using the third cohomology
group.

\vspace{2mm}
 \noindent {\bf Acknowledgement:} We would like to thank M. Markl for his useful comments on free Lie
 $2$-algebras.




\section{Background on Lie $2$-algebras}\label{section2}
\subsection{Basic notions}
$L_{\infty}$-algebras, also called strongly homotopy Lie algebras, were introduced by Drinfeld and Stasheff as a model for ``Lie algebras that satisfy Jacobi identity up to all
higher homotopies''. The following definition of $L_{\infty}$-structure was formulated by Stasheff in 1985. See \cite{Stasheff}.
\begin{defi}\label{defi:L-inf}
An {\bf $L_{\infty}$-algebra} is a graded vector space $\mathfrak{g}=\g_0\oplus{\g_1}\oplus{\cdots}$ equipped with a system $\{l_{k}\mid 1\leq{k}<{\infty}\}$ of linear maps $l_k:\wedge^{k}\g \longrightarrow{\g}$ with degree $deg(l_k)=k-2,$ where the exterior powers are interpreted in the graded sense and the following relation with Koszul sign``Ksgn'' is satisfied for all $n\geq{0}$:
\begin{eqnarray}
\sum_{i+j=n+1}(-1)^{i(j-1)}\sum_{\sigma}sgn(\sigma)Ksgn(\sigma)l_{j}(l_{i}(x_{\sigma(1)},\cdots,x_{\sigma(i)}),x_{\sigma(i+1)},
\cdots,x_{\sigma(n)})=0,
\end{eqnarray}
where the summation is taken over all $(i,n-i)$-unshuffles with $i\geq{1}.$
\end{defi}
Usually, an $n$-term $L_{\infty}$-algebra (i.e., $l_i=0, i\geq n+2$)
is called a Lie $n$-algebra. In particular, if $l_{n+1}=0$, it is
called a strict Lie $n$-algebra. Next we focus on the case of $n=2$.
\begin{defi}\label{defi:Lie 2 homo}
Let $(\g,d,l_{2},l_{3})$ and $(\g{'},d{'},l'_2,l'_3)$ be Lie $2$-algebras. A {\bf Lie $2$-algebra homomorphism} $\varphi:\g \longrightarrow{\g'}$ consists of
\begin{enumerate}
\item[$\bullet$] two linear maps $\varphi_{0}:\g_{0}\longrightarrow{\g'_{0}}$ and $\varphi_{1}:\g_{1}\longrightarrow{\g'_{1}}$,
\item[$\bullet$] one bilinear map $\varphi_{2}: \g_{0}\wedge{\g_{0}}\longrightarrow{\g'_{1}}$,
\end{enumerate}
such that the following equalities hold for all $x,y,z\in{\g_{0}}, a\in{\g_{1}}$ :
\begin{enumerate}
\item[$\bullet$] $d'\circ{\varphi_1}=\varphi_{0}\circ{d}$,
\item[$\bullet$] $\varphi_{0}l_{2}(x,y)-l'_{2}(\varphi_{0}(x),\varphi_{0}(y))=d'\varphi_{2}(x,y)$,
\item[$\bullet$] $\varphi_{1}l_{2}(x,a)-l'_{2}(\varphi_{0}(x),\varphi_{1}(a))=\varphi_{2}(x,da)$,
\item[$\bullet$] $l'_{2}(\varphi_{0}(x),\varphi_{2}(y,z))+c.p.+l'_{3}(\varphi_{0}(x),\varphi_{0}(y),\varphi_{0}(z))=\varphi_{2}(l_{2}(x,y),z)+c.p.+\varphi_{1}(l_{3}(x,y,z))$,
\end{enumerate}
where $c.p.$ means cyclic permutation. It is called a {\bf strong
homomorphism} if $\varphi_{2}=0$.
\end{defi}
\begin{lem}
Let $(\mathfrak{g},d,l_2,l_3)$ be a Lie $2$-algebra and
$\mathfrak{h} \subset \g $ a $2$-vector subspace. Then
$\mathfrak{g}/\mathfrak{h}$ is a quotient Lie $2$-algebra if and
only if
\begin{equation}\label{defi:ideal}
l_2(\mathfrak{h}\wedge
\mathfrak{g})\subset{\mathfrak{h}},\ \ \ \ \ \ \ l_3(\mathfrak{h}_0
\wedge \mathfrak{g}_0\wedge
\mathfrak{g}_0)\subset{\mathfrak{h}_1}.
\end{equation}
\end{lem}
We call $\mathfrak{h}$ satisfying condition (\ref{defi:ideal}) an {\bf ideal} of $\g$.
In fact, the projection
$\pi:\g\longrightarrow\mathfrak{g}/\mathfrak{h}$ becomes a strong
homomorphism. We now give an analogue of the fundamental theorem of
algebras for Lie $2$-algebras.
\begin{pro}\label{pro:ft}
Let $\varphi:\mathfrak{g}\longrightarrow{\mathfrak{g}'}$ be a Lie $2$-algebra homomorphism. Then,
\begin{itemize}
\item[\rm(1)] $\Img\varphi=\Img\varphi_0\oplus\Img\varphi_1$ is a Lie $2$-subalgebra of $\mathfrak{g}'$ if ~$\Img\varphi_2\subset{\Img\varphi_1}$;
\item[\rm(2)] $\ker\varphi=\ker\varphi_0\oplus\ker\varphi_1$ is an ideal of $\mathfrak{g}$ if $\varphi_{2}(\ker \varphi_{0}\wedge \g_0)=0.$
\end{itemize}
Moreover, the two Lie $2$-algebras $\mathfrak{g}/{\ker\varphi}$ and
$\Img\varphi$ are isomorphic if the two conditions above are
satisfied.
\end{pro}
\pf 
By Definition \ref{defi:Lie 2 homo} and $\Img\varphi_2\subset{\Img\varphi_1}$, it is direct to see that $\Img\varphi$ is a 2-vector subspace of $\mathfrak{g}'$ such that $l'_2, l'_3$ are closed on it. Namely, $\Img\varphi$ is a Lie $2$-subalgebra of $\mathfrak{g}'$.

Similarly, by the first three conditions of a Lie 2-algebra homomorphism and $\varphi_{2}(\ker \varphi_{0}\wedge \g_0)=0$, we get $l_{2}(\ker\varphi_0\wedge \mathfrak{g})\subset{\ker\varphi}$. Coupled with the last condition, we further obtain that $l_3(\ker\varphi_0\wedge \mathfrak{g}_0\wedge \mathfrak{g}_0)\subset{\ker\varphi_1}$. That is, $\ker\varphi$ is an ideal of $\mathfrak{g}$. The remaining result is immediate.\qed\vspace{3mm}
\subsection{Cohomology}
 Given a $\g$-module $\mathbb{V}:
V_{1}\stackrel{\partial}{\longrightarrow}V_{0}$ with a Lie 2-algebra
homomorphism $\phi:\g\rightarrow \End(\mathbb{V})$, $\phi$ is called
an {\bf action} of $\g$ on $\mathbb{V}$ and denoted by:
 $$x\triangleright u=\phi(x)u,
 \ \ (x,y)\triangleright u=\phi_2(x,y)u,\ \forall x,y\in \g, u\in \mathbb{V}.$$
The cohomology group of a Lie 2-algebra $(\g,d,l_2,l_3)$ comes from
the generalized Chevalley-Eilenberg complex as follows:
{\footnotesize
\begin{equation} \label{CE complex}
\begin{split}
 &\degree \ -1: V_{1}\stackrel{\D}{\longrightarrow}\\
 &\degree \ \ \ 0\ \ : V_{0}\oplus \Hom(\frkg_{0},V_{1})\stackrel{\D}{\longrightarrow}\\
  &\degree\ \ \ 1\ \ : \Hom(\frkg_{0},V_{0})\oplus\Hom(\frkg_{1},V_{1})\oplus\Hom(\wedge^2\frkg_0,V_{1})\stackrel{\D}{\longrightarrow}\\
  &\degree\ \ \ 2\ \ : \Hom(\frkg_{1},V_{0})\oplus \Hom(\wedge^2\frkg_0,V_{0})\oplus \Hom(\frkg_0\wedge \frkg_{1},V_{1}) \oplus \Hom(\wedge^3\frkg_0,V_{1})\stackrel{\D}{\longrightarrow}\\
  & \stackrel{\D}{\longrightarrow}\cdots.
\end{split}
\end{equation}}
Denote by $C^i(\g,\mathbb{V})$ the set of i-cochains. The coboundary
operator $\D$ can be decomposed as:
\begin{equation*}
  \D=\hat{d}+\hat{\partial}+d^{(1,0)}_\phi+d^{(0,1)}_\phi+d_{\phi_2}+d_{l_3},
\end{equation*}
where, for $s=0,1$,
\begin{eqnarray*}
\hat{d}:\Hom(\wedge^p\frkg_0\wedge \odot
^q\frkg_{1},V_s)&\longrightarrow&\Hom(\wedge  ^{p-1}\frkg_0\wedge
\odot  ^{q+1}\frkg_{1}),V_s),\\
\hat{\partial}:\Hom(\wedge^ p\frkg_0\wedge \odot
^q\frkg_{1},V_{1})&\longrightarrow&\Hom(\wedge^p\frkg_0\wedge
\odot ^q\frkg_{1},V_{0}),\\
d_\phi^{(1,0)}:\Hom(\wedge  ^p\frkg_0\wedge \odot
^q\frkg_{1},V_s)&\longrightarrow& \Hom(\wedge  ^{p+1}\frkg_0\wedge
\odot  ^q\frkg_{1},V_s),\\
d_\phi^{(0,1)}:\Hom(\wedge  ^p\frkg_0\wedge \odot  ^q\frkg_{1},V_0)&\longrightarrow&
\Hom(\wedge  ^p\frkg_0\wedge \odot  ^{q+1}\frkg_{1},V_{1}),\\
d_{\phi_2}:\Hom(\wedge^{p}\frkg_0\wedge \odot^q\frkg_{1},V_0)&\longrightarrow& \Hom(\wedge^{p+2}\frkg_0\wedge \odot^q\frkg_{1},V_{1}),\\
d_{l_3}:\Hom(\wedge^p\frkg_0\wedge \odot^q \frkg_{1},V_s)&\longrightarrow& \Hom(\wedge^{p+3}\frkg_0\wedge \odot^{q-1} \frkg_{1},V_s).
\end{eqnarray*}
More concretely, for any ${x_i}\in{\mathfrak{g}_0},a_i\in{\mathfrak{g}_1}, i\in{\mathbb{N}},$
{\footnotesize
\begin{eqnarray*}
&&\hat{d}f(x_1,\cdots,x_{p-1},a_1,\cdots,a_{q+1})
 =(-1)^{p}\big(f(x_1,\cdots,x_{p-1},dh_1,h_2,\cdots,h_{q+1})+c.p.(a_1,\cdots,a_{q+1})\big),\\
&&\hat{\partial}f=(-1)^{p+2q}\partial\circ f,\\
&&d_\phi^{(1,0)}f(x_1,\cdots,x_{p+1},a_1,\cdots,a_{q})
 =\sum_{i=1}^{p+1}(-1)^{i+1}x_i \triangleright f(x_1,\cdots,\widehat{x_i},\cdots,x_{p+1},
  a_1,\cdots,a_{q})\\
  &&\ \ \ \ \ \ \ \ \ \ \ \ \ \ \ \ \ \ \ \ \ \ \ \ \ \ \ \ \ \ \ \ \ \ \ \ \ \ \ \ \ \ \ \ +\sum_{i<j}(-1)^{i+j}f([x_i,x_j],x_1,\cdots,\widehat{x_i},\cdots,\widehat{x_j},\cdots,x_{p+1},
  a_1,\cdots,a_{q})\\
  &&\ \ \ \ \ \ \ \ \ \ \ \ \ \ \ \ \ \ \ \ \ \ \ \ \ \ \ \ \ \ \ \ \ \ \ \ \ \ \ \ \ \ \ \ +\sum_{i,j}(-1)^{i}f(x_1,\cdots,\widehat{x_i},\cdots,x_{p+1},
  a_1,\cdots,[x_i,a_j],\cdots,a_{q}),\\
&&d_\phi^{(0,1)}f(x_1,\cdots,x_{p},a_1,\cdots,a_{q+1})
 =\sum_{i=1}^{q+1}(-1)^{p}a_i \triangleright f(x_1,\cdots,x_{p},
  a_1,\cdots,\widehat{a_i},\cdots,a_{q+1}),\\
&&d_{\phi_2}f(x_1,\cdots,x_{p+2},a_1,\cdots,a_{q})=
  \sum_{\sigma}(-1)^{p+2q}(-1)^\sigma(x_{\sigma(1)},x_{\sigma(2)})\triangleright f(x_{\sigma(3)},\cdots,x_{\sigma(p+2)},a_1,\cdots,a_{q}),\\
\label{lem:proerty p}&&d_{l_3}f(x_1,\cdots,x_{p+3},a_1,\cdots,a_{q-1})=\sum_\tau-(-1)^\tau f(x_{\tau(4)},\cdots,x_{\tau(p+3)},a_1,\cdots,a_{q-1},l_3(x_{\tau(1)},x_{\tau(2)},x_{\tau(3)})).
\end{eqnarray*}}
where $\sigma$ and $\tau$ are taken over all $(2,p)$-unshuffles and $(3,p)$-unshuffles respectively.

By direct calculations, we get the expressions for 1-cocycles and
1-coboundaries.
\begin{lem}\label{1-cochain}
Let $X+l_X \in C^1(\g,\mathbb{V})$, where $X=(X_0,X_1) \in
\Hom(\g_0,V_0)\oplus\Hom(\g_1,V_1)$ and $l_X\in\Hom(\wedge^2
\g_0,V_1)$. Then
 \begin{itemize}
\item[(1)]
 $\D(X+l_X)=0$ if and only if
\begin{eqnarray}
\label{d} X_0\circ d&=&\partial\circ X_1,\\
\label{a}
\partial l_{X}(x,y)&=&X[x,y]+y\triangleright Xx-x\triangleright Xy,\\
\label{b}l_{X}(x,da)&=&X[x,a]+a\triangleright Xx-x\triangleright Xa,\\
\label{c}Xl_{3}(x,y,z)&=&l_{X}(x,[y,z])+x\triangleright
l_{X}(y,z)-(y,z)\triangleright Xx+c.p.(x,y,z).
\end{eqnarray}
\item[(2)] $\exists u+\Theta\in V_0\oplus
\Hom(\g_0,V_1)=C^0(\g,\mathbb{V})$, s.t. $X+l_X=\D(u+\Theta)$ if and
only if
\begin{eqnarray}
\label{1-cob1}X(x+a)&=&x\triangleright u+a\triangleright u-\partial\Theta(x)-\Theta(da),\\
\label{1-cob2}l_X(x,y)&=&(x,y)\triangleright u+x\triangleright
\Theta(y)-y\triangleright \Theta(x)-\Theta([x,y]).
\end{eqnarray}
 \end{itemize}
\end{lem}

 For future references, we specify the 3-coboundaries. For a cochain $\lambda=\sum_{i=0}^{3}\lambda_i\in
C^2(\g,\mathbb{V})$, where
$$\lambda_0\in{\Hom(\g_{1},V_{0})},
\lambda_1\in{\Hom(\wedge^{2}\g_{0},V_{0})},
\lambda_2\in{\Hom(\g_0\wedge{\g_1},V_{1})}, \lambda_{3}\in
\Hom(\wedge^{3}\g_{0},V_{1}),$$
then $\theta = \D\lambda$
has five components as follows:
  \begin{equation*}
\left\{\begin{array}{rcll} \theta_0&=&d^{(1,0)}_{\phi}\lambda_0+\hat{d}\lambda_1+\hat{\partial}\lambda_2 &~\in{\Hom(\g_{0}\wedge{\g_{1}},V_{0})},\\
\theta_1&=&d^{(0,1)}_{\phi}\lambda_0+\hat{d}\lambda_2&~ \in{\Hom(\odot^2\g_{1},V_{1})},\\
\theta_2&=&d_{l_3}\lambda_0+d^{(1,0)}_{\phi}\lambda_1+\hat{\partial}\lambda_3 &~ \in{\Hom(\wedge^{3}{\g_{0}},V_{0})},\\
\theta_3&=&d_{\phi_{2}}\lambda_0+d^{(0,1)}_{\phi}\lambda_1+d^{(1,0)}_{\phi}\lambda_2+\hat{d}\lambda_3 &~\in{\Hom(\wedge^2{\g_{0}}\wedge{\g_{1}},V_{1})},\\
\theta_4&=& d_{\phi_{2}}\lambda_1+d_{l_3}\lambda_2+d^{(1,0)}_{\phi}\lambda_3 &~ \in{\Hom(\wedge^{4}{\g_{0}},V_{1})}.\end{array}\right.
  \end{equation*}
More precisely, for any $x,y,z,x_i\in \g_0, a,b\in \g_1$,
 \begin{equation}\label{3-ex cochain}
  \left\{\begin{array}{rcll}
\theta_0(x,a)&=&x\triangleright{\lambda_0(a)}-\lambda_0[x,a]+\lambda_1(x,da)-\partial\lambda_2(x,a),\\
\theta_1(a,b)&=&a\triangleright \lambda_0(b)+b\triangleright \lambda_0(a)-\lambda_2(da,b)-\lambda_2(db,a),\\
\theta_2(x,y,z)&=&-\lambda_0l_3(x,y,z)+\big(x\triangleright{\lambda_1(y,z)}-\lambda_1([x,y],z)+c.p.\big)-\partial\lambda_3(x,y,z),\\
\theta_3(x,y,a)&=&(x,y)\triangleright \lambda_0(a)+a\triangleright
\lambda_1(x,y)+x\triangleright{\lambda_2(y,u)}-y\triangleright{\lambda_2(x,a)}-\lambda_2([x,y],a)\\ &&
-\lambda_2(y,[x,a])+\lambda_2(x,[y,a])-\lambda_3(x,y,da),\\
\theta_4(x_1,\cdot \cdot \cdot
,x_4)&=&\sum_{\sigma}(-1)^{\sigma}(x_{\sigma_1},x_{\sigma_2})\triangleright
\lambda_1(x_{\sigma_3},x_{\sigma_4})
-\sum_{\tau}(-1)^{\tau}\lambda_2(x_{\tau_4},l_3(x_{\tau_1},x_{\tau_2},x_{\tau_3}))\\
&& +\sum_{i=1}^{4}(-1)^{i+1}x_i \triangleright
\lambda_3(x_1,\cdots,\widehat{x_i},\cdots,x_{4}) \\ &&
+\sum_{i<j}(-1)^{i+j}\lambda_3([x_i,x_j],x_1,\cdots,\widehat{x_i},\cdots,\widehat{x_j},\cdots,x_{4}).
\end{array}\right.
\end{equation}
\subsection{Derivations}
For a Lie $2$-algebra $(\g,d,[\cdot,\cdot],l_3),$ there is a natural {\bf adjoint action} $ad$ of $\g$ on itself given by
$$ad(x)=[x,\cdot],\ \ \ \ \ ad_2(y,z)=-l_3(y,z,\cdot),\ \ \  \forall x,y,z\in \g.$$

To propose the crossed module of Lie $2$-algebras below, we review
the notion of derivations   of a Lie $2$-algebra $\Der(\g)$, which
was proved to be a strict Lie $2$-algebra in \cite{Sheng1}.

Let $(\g,d,[\cdot,\cdot],l_3)$ be a Lie $2$-algebra. A derivation of
degree 0 of $\mathfrak{g}$ is a pair $(X,l_{X})$, also denoted by
$X+l_X$, where $X=(X_0,X_1)\in\Hom(\g_0,\g_0)\oplus\Hom(\g_1,\g_1)$
and $l_{X}:\g_0 \wedge{\g_0}\rightarrow{\g_1}$ is a linear map, such
that for any $x,y,z\in\g_0,a\in\g_1$,
\begin{equation}\label{neweq}
\left\{\begin{array}{rcll}d\circ X_1&=&X_0\circ d,\\  dl_{X}(x,y)&=&X[x,y]-[Xx,y]-[x,Xy],\\
l_{X}(x,da)&=&X[x,a]-[Xx,a]-[x,Xa],\\
Xl_{3}(x,y,z)&=&l_{X}(x,[y,z])+[x,l_{X}(y,z)]+l_{3}(Xx,y,z)+c.p.(x,y,z).\end{array}\right.
\end{equation}

Denote by $\Der_{0}(\g)$ the set of derivations of degree $0$ of
$\g$. Then one can define a $2$-vector space as
\begin{equation*}
\CD
   \Der(\g): \Der_{1}(\g)\triangleq\Hom(\g_0,\g_1) @>\bar{d}>>\Der_{0}(\g),
\endCD
\end{equation*}
where $\bar{d}$ is given by
$\bar{d}(\Theta)=\delta(\Theta)+l_{\delta(\Theta)}$, in which
$\delta(\Theta)=d\circ \Theta+\Theta\circ d$ and
\[l_{\delta(\Theta)}(x,y)={\Theta}[x,y]-[x,\Theta{y}]-[{\Theta}x,y].\]
In addition, define $\{X+l_{X},\Theta\}=[X,\Theta]$ and
\begin{equation}\label{bra of der}
\{X+l_X,Y+l_Y\}=[X,Y]+X\triangleright l_Y-Y\triangleright l_X,
\end{equation}
where $[\cdot,\cdot]$ is the commutator bracket and
$$X\triangleright l_Y(x,y)=Xl_Y(x,y)-l_Y(Xx,y)-l_Y(x,Xy).$$
\begin{thm}\cite{Sheng1}\label{thm:Der(g)}
With notations above, $(\Der(\g),\{\cdot,\cdot\})$ is a strict Lie
$2$-algebra.
\end{thm}

\begin{rmk}\label{homo}
From the homological viewpoint, we discover an alternative
description of $\Der_0(\g)$ and $\bar{d}$. Comparing Equations
(\ref{neweq}) and $\bar{d}$ with (\ref{d})-(\ref{c}) and
(\ref{1-cob2}) respectively, we note that $\Der_0(\g)$ is indeed the
set of $1$-cocycles of the Lie $2$-algebra $\g$ with respect to the
adjoint action and $\bar{d}=-\D$, which $\D$ is the Lie $2$-algebra
coboundary operator.
\end{rmk}

The adjoint action $ad$ can be extended to a new Lie 2-algebra homomorphism $\overline{ad}:\g\longrightarrow \Der(\g)$, where
\begin{equation}\label{adjoint homo}
\overline{ad}_0(x)=-\D(x)=ad_0(x)+l_3(x,\cdot,\cdot),\ \ \overline{ad}_1=ad_1,
\ \ \overline{ad}_2=ad_2, \ \ \forall x\in\g_0.
\end{equation}

We conclude this section by exploring the derivations of a skeleton
Lie $2$-algebra, which turns out to have an explicit homological
description.
\begin{ex}{\rm
Let $\g_0$ be an Lie algebra and $V$ a $\g_0$-module. Given an Lie algebra
3-cocycle $l_3\in C^3(\g_0,V)$, we get a skeletal Lie $2$-algebra
$\g=(V\stackrel{0}\longrightarrow {\g_0},l_2,l_3)$, where $l_2$ is
defined by $$l_2^0(x,y)=[x,y]_{\g_0},\ \ \ \ \
l_2^1(x,u)=x\triangleright u,\ \ \forall x,y\in \g_0,\ u\in V.$$

For $X=(X_0,X_1)\in \Hom(\g_0,\g_0)\oplus\Hom(V,V)$
and $l_X\in C^2(\g_0,V)$,  it is easy to check that
\begin{equation*}
X+l_X\in\Der_0(\g) \Longleftrightarrow \left\{\begin{array}{rcll}X&\in& \Der(\g_0\ltimes V),\\
\mathfrak{D} l_X&=&[X,l_3],\end{array}\right.
\end{equation*}
where $\mathfrak{D}: C^k(\g_0,V)\rightarrow C^{k+1}(\g_0,V)$ is the
Lie algebra coboundary operator,  $\Der(\g_0\ltimes V)$ is the Lie
algebra of derivations of the semi-product Lie algebra $\g_0\ltimes
V$ and the bracket $[\cdot,\cdot]$ is given by
$$[X,l_3](x,y,z)=X_1l_3(x,y,z)-l_3(X_0x,y,z)-l_3(x,X_0y,z)-l_3(x,y,X_0z),\ \ \ \forall x,y,z\in \g_0.$$
In fact, such a bracket was used to introduce the notion of pre-Lie
algebras by Gerstenhaber in \cite{Gerstenhaber2}. In this case, we
have $\Der_1(\g)=C^1(\g_0,V)$ and the map $\bar{d}:
\Der_1(\g)\rightarrow \Der_0(\g)$ is given by
$\bar{d}(\Theta)=0-\mathfrak{D}(\Theta)$.}
\end{ex}
\section{Crossed modules of Lie $2$-algebras}\label{section cm}
\subsection{Definition of crossed modules}
 Let $(\mathfrak{m}:\mathfrak{m}_{1}\stackrel{\tilde{d}}{\longrightarrow}\mathfrak{m}_{0},\tilde{l}_2,\tilde{l}_3)$ and $(\mathfrak{g}:\mathfrak{g}_{1}\stackrel{d}{\longrightarrow}\mathfrak{g}_{0},l_2,l_3)$ be two Lie 2-algebras.
 We call {\bf $\g$ acts on $\m$ by derivations} if there exists a Lie $2$-algebra homomorphism $\phi:\g\longrightarrow \End(\m)$ and a linear map $l_{\phi_0(x)}:\wedge^2\m_0\longrightarrow\m_1$ such that $\phi_0(x)+l_{\phi_0(x)}\in \Der_0(\m)$ and the map
 $$(\phi_0+l_{\phi_0},\phi_1,\phi_2):\g\longrightarrow \Der(\m)$$ is a Lie $2$-algebra homomorphism. By abuse of notations, we denote by $\phi$ both the action and the action by derivations. Then we shall define a {\bf crossed product} of $\g$ and $\m$ denoted by $\g\triangleright_\phi \m$, which is still a Lie $2$-algebra depending on the following lemma.

\begin{lem}\label{lem:crossed product}
Let $\phi$ be an action of $\mathfrak{g}$ on $\mathfrak{m}$ by
derivations, then $\g\triangleright_\phi
\m\triangleq(\mathfrak{g}\oplus{\mathfrak{m}},L_1,[\cdot,\cdot],L_3)$
is a Lie $2$-algebra with $\mathfrak{g}$ as a lie $2$-subalgebra and
$\mathfrak{m}$ as an ideal, where $L_1=d+\tilde{d}$ and
\begin{equation*}
\left\{\begin{array}{rcl}
[x+\alpha,y+\beta]&=&{l_2(x,y)}+\tilde{l}_2(\alpha,\beta)+x\triangleright{\beta}-y\triangleright{\alpha},$\ \ \ \ $\forall{x,y}\in{\mathfrak{g}},\forall{\alpha,\beta}\in{\mathfrak{m}},\\
L_3(x+\alpha,y+\beta,z+\gamma)&=&{l_3(x,y,z)}+\tilde{l}_3(\alpha,\beta,\gamma)
    -(x,y)\triangleright \gamma -(y,z)\triangleright \alpha\\  &&
    -(z,x)\triangleright \beta+l_{\phi_{0}(x)}(\beta,\gamma)+l_{\phi_{0}(y)}(\gamma,\alpha)
    +l_{\phi_{0}(z)}(\alpha,\beta),\\ && \forall{x,y,z}\in{\mathfrak{g}_0},\forall{\alpha,\beta,\gamma}\in{\mathfrak{m}_0}.\end{array}\right.
\end{equation*}

Conversely, let $(\theta,L_1,L_2,L_3)$ be a Lie $2$-algebra which
can be split into the direct sum of a Lie $2$-subalgebra
$\mathfrak{g}$ and an ideal $\mathfrak{m}$, then there exists an
action $\phi$ of $\g$ on $\m$ by derivations such that $\theta=\g\triangleright_\phi
\m$, where
$\phi:\mathfrak{g}\longrightarrow{\Der(\mathfrak{m})}$ is defined by
\begin{equation*}
\left\{\begin{array}{rcll} \phi_{0}(x)+l_{\phi_{0}(x)}&=&L_2(x,\cdot)+L_3(x,\cdot,\cdot),&~\forall x\in\g_0,\\
\phi_{1}(a)&=&{L_2(a,\cdot)},& ~\forall a\in\g_1,\\
 \phi_{2}(x,y)&=&-L_3(x,y,\cdot),& ~\forall x,y\in\g_0.\end{array}\right.
\end{equation*}
\end{lem}

\pf We merely provide the proof of the coherence law of $L_2,L_3$.
Denote by $[\cdot,\cdot]$ both $l_2$ and $\tilde{l}_2$ if there is
no risk of confusion. Firstly, for four elements in $\mathfrak{g}_0$
or four elements in $\mathfrak{m}_0$, it holds since $\mathfrak{g}$
and $\mathfrak{m}$ are Lie 2-algebras and $L_2,L_3$ reserve their brackets. For any
$x,y,z\in{\mathfrak{g}_0},\alpha,\beta,\gamma\in{\mathfrak{m}_0}$, we have
\begin{eqnarray*}
&&[L_{3}(x,\alpha,\beta),\gamma]+L_{3}(x,[\alpha,\beta],\gamma)-L_{3}([x,\alpha],\beta,\gamma)+c.p.(\alpha,\beta,\gamma)-[L_{3}(\alpha,\beta,\gamma),x]\\
&=&
[l_{\phi_{0}(x)}(\alpha,\beta),\gamma]+l_{\phi_{0}(x)}([\alpha,\beta],\gamma)-\tilde{l}_3
(x\triangleright{\alpha},\beta,\gamma)+c.p.(\alpha,\beta,\gamma)
+x\triangleright{\tilde{l}_3(\alpha,\beta,\gamma)},
\end{eqnarray*}
which vanishes since $\phi_{0}(x)+l_{\phi_{0}(x)}\in \Der_0(\m)$. Next,
making use of the fact $\phi$ is an action by derivations coupled with (\ref{bra of der}), we have
\[l_{\phi_{0}[x,y]}-(\phi_{0}(x)\triangleright l_{\phi_{0}(y)}-\phi_{0}(y)\triangleright l_{\phi_{0}(x)})=l_{\delta\phi_{2}(x,y)}.\]
Hence,
\begin{eqnarray*}
&&-[L_{3}(y,\alpha,\beta),x]+L_{3}([x,\alpha],y,\beta)-L_{3}([x,\beta],y,\alpha)-c.p.(x,y)\\
&&+
[L_{3}(x,y,\alpha),\beta]-[L_{3}(\beta,x,y),\alpha]-L_{3}([\alpha,\beta],x,y)-L_{3}([x,y],\alpha,\beta)\\
&=&
x\triangleright{l_{\phi_{0}(y)}(\alpha,\beta)}-l_{\phi_{0}(y)}(x\triangleright{\alpha},\beta)-l_{\phi_{0}(y)}(\alpha,x\triangleright{\beta})-c.p.(x,y)
\\ &&-[(x,y)\triangleright \alpha,\beta]+[(x,y)\triangleright \beta,\alpha]+(x,y)\triangleright [\alpha,\beta]-l_{\phi_{0}[x,y]}(\alpha,\beta)
\\ &=&(\phi_{0}(x)\triangleright l_{\phi_{0}(y)}-\phi_{0}(y)\triangleright l_{\phi_{0}(x)})(\alpha,\beta)+l_{\delta\phi_{2}(x,y)}(\alpha,\beta)-l_{\phi_{0}[x,y]}(\alpha,\beta)
\\ &=&0.
\end{eqnarray*}
Now, it remains to show
\begin{eqnarray*}
&&-[L_{3}(y,z,\alpha),x]-L_{3}([x,\alpha],y,z)-L_{3}([x,y],z,\alpha)+c.p.(x,y,z)+[L_{3}(x,y,z),\alpha]
\\&=&-x\triangleright{((y,z)\triangleright\alpha)}+(y,z)\triangleright (x\triangleright{\alpha})+([x,y],z)\triangleright \alpha+c.p.(x,y,z)+[L_{3}(x,y,z),\alpha]
\\&=&\big([\phi_{2}(y,z),\phi_{0}(x)]+\phi_{2}([x,y],z)+c.p.(x,y,z)+\phi_{1}(l_{3}(x,y,z))\big)\alpha
\\&=&0,
\end{eqnarray*}
where the last equality follows from that $\phi$ is a Lie 2-algebra
homomorphism. This finishes the proof of $\g\triangleright_\phi \m$
is Lie 2-algebra. The remaining results are easy to
get.\qed\vspace{3mm}

\begin{defi}\label{defi:Lie 2 cm}
A { \bf crossed module of Lie 2-algebras} is a quadruple
$(\m,\g,\phi,\Pi)$, where $\m,\g$ are two Lie $2$-algebras, $\phi$
is an action of $\g$ on $\m$ by derivations, and
$\Pi:\g\triangleright_\phi \m \rightarrow  \g $ is a Lie $2$-algebra
homomorphism, such that $\Pi|_{\g}=Id=(id,id,0)$ and
\begin{itemize}
\item[\rm(i)] $\tilde{l}_2(\alpha,\beta)=\Pi(\alpha)\triangleright{\beta},\ \ \ \ \forall \alpha,\beta\in{\mathfrak{m}}$,
\item[\rm(ii)] $\tilde{l}_3(\alpha,\beta,\gamma)=-(\Pi_0 \alpha,\Pi_0 \beta)\triangleright \gamma
    -\Pi_2(\Pi_0\alpha,\beta)\triangleright{\gamma},\ \ \ \ \forall{\alpha,\beta,\gamma}\in{\mathfrak{m}_0},$
\item[\rm(iii)] $l_{\phi_0(x)}(\beta,\gamma)=-(x,\Pi_0 \beta)\triangleright \gamma
    -\Pi_2(x,\beta)\triangleright{\gamma},\ \ \ \ \forall{\beta,\gamma}\in{\mathfrak{m}_0},x\in{\mathfrak{g}_0}$,
\item[\rm(iv)]
    $\Pi_2(\alpha,\beta)=\Pi_2(\Pi_{0}\alpha,\beta)=\Pi_2(\alpha,\Pi_{0}\beta)
   ,\ \ \ \ \forall{\alpha,\beta}\in{\mathfrak{m}_0}$.
\end{itemize}
In particular, it is called a {\bf strong crossed module of Lie
$2$-algebras} if ${\Pi_2}=0$.
\end{defi}
We drop the words ``Lie $2$-algebras'' except when emphasis is
needed. For a crossed module $(\m,\g,\phi,\Pi)$, decompose $\Pi$
into
$$\Pi=(\Pi_0,\Pi_1,\Pi_2)=Id + \sigma + \varphi=\big((id,\varphi_0),(id,\varphi_1),(0,\sigma,\varphi_2)\big)$$ where $\varphi=\Pi|_{\m}$ and $\sigma={\Pi_2}|_{\g_0\wedge \m_0}$. It is evident that $\varphi:\m\longrightarrow \g$ is a Lie 2-algebra homomorphism. In the following, we always describe a crossed module as $(\m,\g,\phi,\varphi,\sigma)$. Then a strong crossed module means $\sigma=0$ and is denoted by $(\m,\g,\phi,\varphi).$
\begin{defi}
Let $(\m,\g,\phi,\varphi,\sigma)$ and
$(\m',\g',\phi',\varphi',\sigma')$ be crossed modules, a {\bf
morphism of crossed modules} consists of two Lie $2$-algebra
homomorphisms $F:\m\longrightarrow \m'$, $G:\g\longrightarrow \g'$
and a linear map $\tau:\g_0\wedge \m_0\longrightarrow \m'_1$ such
that $\varphi'\circ F=G\circ\varphi$ and
$$\big((G_0,F_0),(G_1,F_1),(G_2,\tau,F_2)\big):\g\triangleright_{\phi}\m\longrightarrow
\g'\triangleright_{\phi'}\m'$$ is a homomorphism of Lie $2$-algebras.
If $G_2,F_2,\tau$ vanish, we call it a strong morphism.
\end{defi}

Similar to the Lie algebra case, we have the following proposition.
\begin{pro}\label{pro:Lie 2 cm}
Let $\g$ be a Lie $2$-algebra and $\m$ a $\g$-module. Given a chain
map $\varphi:\m\rightarrow \g$ and a map
$\sigma:\g_0\wedge{\m_0}\longrightarrow{\g_1}$, satisfying that
$$\Pi\triangleq
\big((id,\varphi_0),(id,\varphi_1),(0,\sigma,0)\big):\mathfrak{g}\ltimes_\phi{\mathfrak{m}}\longrightarrow
\g$$ is a Lie $2$-algebra homomorphism and
\begin{itemize}
\item[\rm(1)] $\varphi(\alpha)\triangleright{\beta}=-\varphi(\beta)\triangleright{\alpha},\ \ \ \ \forall{\alpha,\beta}\in{\mathfrak{m}}$,
\item[\rm(2)] $(\varphi_{0}\alpha,\varphi_{0}\beta)\triangleright \gamma+\sigma(\varphi_{0}\alpha,\beta)\triangleright
   {\gamma}=
-(\varphi_{0}\alpha,\varphi_{0}\gamma)\triangleright
\beta-\sigma(\varphi_{0}\alpha,\gamma)\triangleright {\beta},\ \ \ \
\forall{\alpha,\beta,\gamma}\in{\mathfrak{m}_0}$,
\item[\rm(3)] $(x,\varphi_{0}\beta)\triangleright \gamma+\sigma(x,\beta)\triangleright{\gamma}=
-(x,\varphi_{0}\gamma)\triangleright
\beta-\sigma(x,\gamma)\triangleright{\beta},\ \ \ \ \
\forall{\beta,\gamma}\in{\mathfrak{m}_0},x \in{\mathfrak{g}_0}$,
\item[\rm(4)] $\sigma(\varphi_{0}\alpha,\beta)=\sigma(\alpha,\varphi_{0}\beta),\ \ \ \ \ \forall{\alpha,\beta}\in{\mathfrak{m}_0},$
\end{itemize}
then, there exists a unique Lie $2$-algebra structure on
$\mathfrak{m}$, linear maps $l_{\phi_0(x)}:\wedge^{2}\m_0
\longrightarrow{\m_1}$ and $\varphi_2:\wedge^{2}\m_0
\longrightarrow{\mathfrak{g}_1}$ such that
$(\m,\g,\hat{\phi},\hat{\varphi},\sigma)$ is a crossed module, where
$\hat{\phi}=(\phi_0+l_{\phi_0},\phi_1,\phi_2)$ and
$\hat{\varphi}=(\varphi_0,\varphi_1,\varphi_2)$.
\end{pro}
\pf Define $\tilde{l}_2,\tilde{l}_3,l_{\phi_0(x)},\varphi_2$ on $\m$
by the right hand sides of equalities $(i)$-$(iv)$ of Definition
\ref{defi:Lie 2 cm}. Then, by direct verification, we obtain that
$\mathfrak{m}$ with $\tilde{l}_2,\tilde{l}_3$ is a Lie 2-algebra and
$\hat{\phi}:\g\longrightarrow \Der(\m)$ is a Lie 2-algebra
homomorphism. Moreover,
$\hat{\Pi}=((id,\varphi_0),(id,\varphi_1),(0,\sigma,\varphi_2)):
\mathfrak{g}\triangleright_{\hat{\phi}}{\mathfrak{m}}\longrightarrow
\g$ is a Lie 2-algebra homomorphism. Thus, we get a crossed module
of Lie 2-algebras.\qed\vspace{3mm}
\subsection{Lie 3-algebras associated to crossed modules}
A crossed module $(\m,\g,\phi,\varphi,\sigma)$ corresponds to a square
\begin{equation*}
\xymatrix{ \m_1 \ar@{->} [r]^{\varphi_1}\ar@ {->} [d]_{\tilde{d}}
&\g_1 \ar@ {->} [d]^{d}\\
{\m_0} \ar@{->} [r]_{\varphi_0}&\g_0. }
\end{equation*}
Consider its mapping cone complex (\cite{Loday})
\begin{equation}
\xymatrix@C=0.5cm{\mathfrak{V} : \mathfrak{m}_{1}\ar[rr]^{d_D} &&
{\mathfrak{g}_{1}\oplus{\mathfrak{m}_{0}}} \ar[rr]^{d_D} &&
                \mathfrak{g}_{0},}
\end{equation}
where
$$d_D(\xi)=-\varphi_{1}\xi+\tilde{d}\xi, \ \ \ d_D(a+\alpha)=da+\varphi_{0}\alpha, \, ~~~  \, \forall\xi\in{\mathfrak{m}_{1}}, a+\alpha
\in{\mathfrak{g}_{1}\oplus{\mathfrak{m}_0}}.$$ Then $d_D^2=0$
follows from the fact that
$d\circ{\varphi_{1}}=\varphi_{0}\circ{\tilde{d}}.$

Define $\Courant{\cdot,\cdot}$ and $l^3$ on it by: for any $x,y,z\in
\g_0$, $a,b\in \g_1$, $\alpha,\beta\in \m_0$ and $\xi\in \m_1$,
\begin{equation*}
\left\{\begin{array}{rcl} \Courant{x,y}&=&l_{2}(x,y),\\ \Courant{x,\xi}&=&-\Courant{\xi,x}=x\triangleright{\xi},\\ \Courant{a+\alpha,b+\beta}&=&a\triangleright{\beta}+b\triangleright{\alpha},\\
\Courant{x,a+\alpha}&=&-\Courant{a+\alpha,x}=l_2(x,a)-\sigma(x,\alpha)+x\triangleright{\alpha}, \\
l^3(x,y,z)&=&l_{3}(x,y,z),\\
l^{3}(x,y,a+\alpha)&=&-l^{3}(x,a+\alpha,y)=l^{3}(a+\alpha,x,y)=-(x,y)\triangleright
\alpha.\end{array}\right.
\end{equation*}
\begin{thm}\label{thm:Lie 3}
With the above notations, $(\mathfrak{V}, d_D,
\Courant{\cdot,\cdot}, l^3)$ is a strict Lie $3$-algebra.
\end{thm}
\pf We need to verify all the conditions of Lie 3-algebras. Firstly, it is obvious that
$\Courant{\cdot,\cdot}$ and $l^3$ are antisymmetric in the graded
sense. Then, it suffices to prove that equality $(1)$ of
Definition \ref{defi:L-inf} holds for $1\leq n\leq 5.$

$\bullet$ $n=1$: $(1)$ reduces to $d_D^{2}=0$, which has already
been checked.

$\bullet$ $n=2$: Condition $(1)$ gives
$$d_D[x,y]=\Courant{d_Dx,y}+(-1)^{|x|}\Courant{x,d_Dy},\ \ \ \ \ \forall x,y \in \mathfrak{V},$$ which is equivalent to
\begin{equation*}
\left\{\begin{array}{rcl} d_{D}\Courant{x,\xi}&=&\Courant{x,d_{D}\xi},\\
d_{D}\Courant{x,a+\alpha}&=&\Courant{x,d_{D}(a+\alpha)},\\
\Courant{d_{D}(a+\alpha),\xi}&=&\Courant{a+\alpha,d_{D}\xi},\\
d_D\Courant{a+\alpha,b+\beta}&=&\Courant{d_{D}(a+\alpha),b+\beta}
-\Courant{a+\alpha,d_{D}(b+\beta)}.\end{array}\right.
\end{equation*}
The first three equations are easy to verify. As for the last one,
by direct computations, we have
\begin{eqnarray*}
d_D\Courant{a+\alpha,b+\beta}&=&d_D(a\triangleright
\beta+b\triangleright \alpha)
\\ &=&-\varphi_{1}(a\triangleright \beta)-\varphi_{1}(b\triangleright \alpha)+\tilde{d}(a\triangleright \beta)+\tilde{d}(b\triangleright \alpha),
\end{eqnarray*}
and
\begin{eqnarray*}
&&\Courant{d_{D}(a+\alpha),b+\beta}-\Courant{a+\alpha,d_{D}(b+\beta)}\\
&=&[da+\varphi_{0}\alpha,b]-\sigma(da+\varphi_{0}\alpha,\beta)+(da+\varphi_{0}\alpha)\triangleright\beta\\
&&
+[db+\varphi_{0}\beta,a]-\sigma(db+\varphi_{0}\beta,\alpha)+(db+\varphi_{0}\beta)\triangleright\alpha\\
&=&[\varphi_{0}\alpha,b]+\sigma(\alpha,db)-[a,\varphi_{0}\beta]-\sigma(da,\beta)+
da\triangleright \beta+db\triangleright \alpha,
\end{eqnarray*}
where we have used conditions $(i)$ and $(iv)$ of Definition
\ref{defi:Lie 2 cm}. Therefore, the equality
$d_D\Courant{a+\alpha,b+\beta}=\Courant{d_{D}(a+\alpha),b+\beta}
-\Courant{a+\alpha,d_{D}(b+\beta)}$ holds since $\Pi$ is a
homomorphism.

$\bullet$ $n=3$: We are supposed to check the graded Jacobi
identity: for any $x,y,z\in \mathfrak{V},$
\begin{eqnarray*}
&&(-1)^{|x|\cdot|z|}\Courant{\Courant{x,y},z}+c.p.\\ &=&
(-1)^{|x|\cdot|z|+1}\big\{d_Dl^{3}(x,y,z)+l^{3}(d_Dx,y,z)+(-1)^{|x|}l^{3}(x,d_Dy,z)
+(-1)^{|x|+|y|}l^{3}(x,y,d_Dz)\big\}.
 \end{eqnarray*}
Following from that $\triangleright$ is an action and $\Pi$ is
a homomorphism, we have
\[\Courant{\Courant{x,y},z}+c.p.=-dl_3(x,y,z)=-d_Dl^3(x,y,z),\]
\[\Courant{\Courant{x,y},\xi}+c.p.=(x,y)\triangleright \tilde{d}\xi=-l^3(x,y,d_{D}\xi),\]
and
\begin{eqnarray*}
&&\Courant{\Courant{x,y},a+\alpha}+c.p.\\
&=&[[x,y],a]-\sigma([x,y],\alpha)+[x,y]\triangleright\alpha+
[[y,a],x]-[\sigma(y,\alpha),x]-\sigma(y\triangleright\alpha,x)-x\triangleright(y\triangleright\alpha)\\
&&
-[[x,a],y]+[\sigma(x,\alpha),y]+\sigma(x\triangleright\alpha,y)+y\triangleright(x\triangleright\alpha)\\
&=&-l_{3}(x,y,da)+\tilde{d}((x,y)\triangleright
\alpha)-\varphi_{1}((x,y)\triangleright
\alpha)-l_3(x,y,\varphi_{0}\alpha)
\\ &=&-l^3(x,y,d_{D}(a+\alpha))-d_{D}l^3(x,y,a+\alpha).
\end{eqnarray*}
 The next case is,
\begin{eqnarray*}
&&\Courant{\Courant{x,a+\alpha},b+\beta}+\Courant{\Courant{a+\alpha,b+\beta},x}-\Courant{\Courant{b+\beta,x},a+\alpha}
\\ &=&[x,a]\triangleright\beta-\sigma(x,\alpha)\triangleright\beta+b\triangleright(x\triangleright\alpha)
-x\triangleright(a\triangleright\beta)-x\triangleright(b\triangleright\alpha)\\
&&+[x,b]\triangleright\alpha-\sigma(x,\beta)\triangleright\alpha+a\triangleright(x\triangleright\beta)\\
&=&(x,da)\triangleright \beta+(x,db)\triangleright \alpha
+(x,\varphi_{0}\alpha)\triangleright \beta+(x,\varphi_{0}\beta)\triangleright \alpha\\
&=&-l^3(x,d_{D}(a+\alpha),b+\beta)+l^3(x,a+\alpha,d_{D}(b+\beta)),
\end{eqnarray*}
where we have used the equation $\sigma(x,\alpha)\triangleright
\beta+\sigma(x,\beta)\triangleright \alpha
=-(x,\varphi_{0}\alpha)\triangleright
\beta-(x,\varphi_{0}\beta)\triangleright \alpha$ followed from
condition $(iii)$ of Definition \ref{defi:Lie 2 cm}. This finishes
the proof of the graded Jacobi identity.

$\bullet$ $n=4$: 
Specifically, for four elements in $\mathfrak{V}_0$, $(1)$ holds
since $\mathfrak{g}$ is a Lie 2-algebra and the definition of $l^3$.
While for three elements in $\mathfrak{V}_0$ and one element in
$\mathfrak{V}_1$, by straightforward deduce, condition $(1)$ is
equivalent to the coherence law of $L_2$ and $L_3$ on three elements
in $\mathfrak{g}_0$ and one element in $\mathfrak{m}_0$ in Lemma
\ref{lem:crossed product}. By careful analysis, all the other cases
are trivial.

$\bullet$ $n=5$: We shall prove
$$\Sigma_{\sigma}(-1)^{\sigma}Ksgn(\sigma)l^3(l^3(x_{\sigma_1},x_{\sigma_2},x_{\sigma_3}),x_{\sigma_4},x_{\sigma_5})=0.$$
Actually, every term in the summation vanishes by the definition of
$l^3$.
This completes the proof.\qed\vspace{3mm}
\subsection{Examples}
\begin{ex}\label{ex:Lie Lie 2 cm}{\rm

Let $(\m,\g,\phi,\varphi)$ be a strong crossed module with strict
Lie $2$-algebras $\m,\g$ and strong homomorphism $\phi$. Treating
$\g$ as a crossed module of Lie algebras with $[a,b]=[da,b]$ on
$\g_1$ and $\m$ likewise, we get a commutative diagram
\begin{equation*}
\xymatrix{ \m_1 \ar@{->} [r]^{\varphi_1}\ar@ {->} [d]_{\tilde{d}}
\ar@{->}[dr]^{\varphi_0\circ\tilde{d}}
&\g_1 \ar@ {->} [d]^{d}\\
{\m_0} \ar@{->} [r]_{\varphi_0}&\g_0 }
\end{equation*}
such that all the maps are crossed modules of Lie algebras.
Moreover, defining $[a,\xi]=da\triangleright\xi,
\forall a\in \g_1,\xi\in \m_1$, the crossed product $\g\triangleright_\phi \m$ is also a
crossed module of Lie algebras.}
\end{ex}
\begin{rmk}
This example remains us of the notion of crossed squares and $2$-crossed modules of Lie algebras introduced by Ellis in \cite{Elli}. The group-theoretic setting is due to Conduch${\rm\acute{e}}$ (\cite{Condu}). See also \cite{Mu, Faria} for more details. The relation between the strong crossed modules in the above example and crossed squares of Lie algebras is still a mystery to us, which deserves to be further studied.
\end{rmk}
\begin{ex} \label{ex1} {\rm For a Lie $2$-algebra $(\mathfrak{g},d,l_{2},l_{3})$, it is obvious that $\Der(\mathfrak{g})$ acts on $\mathfrak{g}$ by derivations with $Id:\Der(\g)\longrightarrow\Der(\g)$. Consider the adjoint homomorphism $\overline{ad}:\mathfrak{g}\longrightarrow \Der(\mathfrak{g})$ given by (\ref{adjoint homo}) and a linear map $\sigma:\Der_0(\g)\wedge\g_0\longrightarrow \Der_1(\g)$ defined by
\[\sigma(X+l_X,x)=-l_X(x,\cdot), \ \ \ \ \ \forall X+l_X\in{\Der_0(\mathfrak{g})}, {x}\in{\mathfrak{g}_{0}}.\]
By straightforward verification, we get
$(\g,\Der(\g),Id,\overline{ad},\sigma)$ is a crossed module. Note
that this crossed module is not strong even if $\g$ is a strict Lie
$2$-algebra. From Theorem \ref{thm:Lie 3} it follows that the
$3$-term complex of vector spaces
\begin{equation*}
\xymatrix@C=0.5cm{\DER(\mathfrak{g}) : \mathfrak{g}_1\ar[rr]^{d_{D}}
&& {\Der_1(\mathfrak{g})\oplus{\mathfrak{g}_0}} \ar[rr]^{d_{D}} &&
                \Der_0(\mathfrak{g})
                }
\end{equation*}
is a strict Lie $3$-algebra. Moreover, $l^3$ vanishes by definition. This recovers \cite[Theorem 3.8]{Sheng1}.}
\end{ex}
\begin{ex}\label{ideal Lie 2 cm}{\rm
Let $(\mathfrak{g},d,l_2,l_3)$ be a Lie $2$-algebra and
$\mathfrak{m}$ an ideal of $\g$, then $(\m,\g,\overline{ad},i)$ is a
strong crossed module, in which $\overline{ad}:\g\longrightarrow
\Der(\m)$ given by (\ref{adjoint homo}) is an action of
$\mathfrak{g}$ on $\mathfrak{m}$ by derivations and the inclusion
map $i:\m\longrightarrow \g$ is a strong homomorphism.}
\end{ex}
According to Example \ref{ex1} and Example \ref{ideal Lie 2 cm} with
$\m=\g$, there are two natural crossed modules for a Lie $2$-algebra
$\g$. Moreover, we get a commutative diagram
\begin{equation*}
\CD
   \g\triangleright_{\overline{ad}}\g @>\overline{ad}\oplus Id>> \Der(\g)\triangleright_{Id}\g \\
  @V Id+IdVV @VVId+\sigma+\overline{ad} V \\
  \g @>>\overline{ad}> \Der(\g)
\endCD
\end{equation*}
such that all the maps are Lie $2$-algebra homomorphisms. Namely,
$(Id,\overline{ad},0)$ is a homomorphism between the two crossed
modules.

\begin{ex}\label{2-vector spaces}{\rm The following example is inspired by \cite[Example $3$]{Wagemann}.
Given a Lie $2$-algebra $\mathfrak{h}$, a short exact sequence of
$\mathfrak{h}$-modules:
\begin{equation}\label{ex se1}
0\longrightarrow \mathbb{V} \stackrel{p}\longrightarrow
\mathbb{I}\stackrel{q} \longrightarrow \mathbb{Q}\longrightarrow0,
\end{equation}
(regarded as a short exact sequence of trivial Lie $2$-algebras) and
a $2$-cocycle $\lambda\in C^2(\mathfrak{h},\mathbb{Q})$, by
\cite[Theorem 4.5]{Liu}, we get an abelian extension of
$\mathfrak{h}$ by $\mathbb{Q}$:
\begin{equation}\label{ex se2}
 0\rightarrow \mathbb{Q}
\stackrel{i}\longrightarrow \mathfrak{h}\oplus_\lambda \mathbb{Q}
\stackrel{\pi} \longrightarrow \mathfrak{h}\rightarrow 0.
\end{equation}
 Splicing (\ref{ex se1}) and (\ref{ex se2}) together, we get a strong crossed module $\varepsilon=(\mathbb{I},\mathfrak{h}\oplus_\lambda \mathbb{Q},\phi,\varphi)$,
 where $\mathbb{I}\stackrel{\varphi}\longrightarrow \mathfrak{h}\oplus_\lambda \mathbb{Q}$ is defined by $\varphi(v)=(0, q(v))$,
 and $\phi$ is an action of $\mathfrak{h}\oplus_\lambda \mathbb{Q}$ on $\mathbb{I}$ given by the action
 of $\mathfrak{h}$ on $\mathbb{I}$. Since $\mathbb{I}$ is a trivial Lie $2$-algebra, by defining $l_{\phi_0}=0$, it is clear that $\phi$ is an action by derivations.}
\end{ex}

\section{The first cohomology and derivations of crossed modules}\label{1-degree}
\subsection{Lie algebra structures on $C^1(\g,\mathbb{V})$}
Let $\mathbb{V}:V_{1}\stackrel{\partial}{\longrightarrow}{V_{0}}$ and $ \mathbb{W}:W_{1}\stackrel{d}{\longrightarrow}{W_{0}}$ be two 2-vector spaces. Then we can construct a new 2-vector space
$$\Hom(\mathbb{V},\mathbb{W}):\Hom_1(\mathbb{V},\mathbb{W})\stackrel{\delta}\longrightarrow \Hom_0(\mathbb{V},\mathbb{W})$$
where $\Hom_1(\mathbb{V},\mathbb{W})=\Hom(V_0,W_1)$,
$$ \Hom_0(\mathbb{V},\mathbb{W})=\{X_0+X_1\in \Hom(V_0,W_0)\oplus \Hom(V_1,W_1)|X_0\circ \partial= d\circ X_1\},$$
and $\delta(\Theta)= d\circ \Theta+\Theta\circ
\partial$.
\begin{lem}\label{lem:END}
Given $\varphi\in \Hom_0(\mathbb{W},\mathbb{V})$, there exists a
strict Lie $2$-algebra structure on $\Hom(\mathbb{V},\mathbb{W})$,
where the bracket $[\cdot,\cdot]_\varphi$ is given by:
\begin{eqnarray}\label{eq:phi bracket}
[X,Y]_\varphi=X\circ \varphi\circ Y-Y\circ \varphi\circ X, \ \ \
[X,\Theta]_\varphi=X\circ \varphi\circ \Theta-\Theta\circ
\varphi\circ X,\ \ \ [\Theta,\Theta']_\varphi=0.
\end{eqnarray}
\end{lem}

In particular, if $\mathbb{W}=\mathbb{V},\varphi=Id$, this recovers the strict Lie 2-algebra $\End(\mathbb{V})$, which plays the same role as $\gl(V)$ for a vector space $V$.

Set
$\La(\g,\mathbb{V})\triangleq
\Hom(\frkg_{0},V_{0})\oplus\Hom(\frkg_{1},V_{1})$, then $C^1(\g,\mathbb{V})=\La(\g,\mathbb{V})\oplus\Hom(\wedge^2\frkg_0,V_{1})$. Given an element $\varphi=(\varphi_0,\varphi_1)\in \Hom(V_{0},\g_0)\oplus\Hom(V_{1},\g_1)$, we have $\La(\g,\mathbb{V})$ is a Lie algebra with the bracket $[\cdot,\cdot]_\varphi$ given by the first formula of (\ref{eq:phi bracket}). Furthermore, we have:
\begin{lem}\label{lem:degree 1 Lie alg}
\begin{itemize}
\item[\rm(1)] $\Hom(\wedge^2\frkg_0,V_{1})$ is an $\La(\g,\mathbb{V})$-module, where the action is given by
\begin{eqnarray*}
(X\triangleright \xi)(x,y)=X_1\varphi_1 \xi(x,y)-\xi(\varphi_0 X_0 x,y)-\xi(x,\varphi_0 X_0y),\ \ \forall x,y\in{\g_0},
\end{eqnarray*}
for any $X=(X_0,X_1)\in \La(\g,\mathbb{V})$ and $\xi\in \Hom(\wedge^2\frkg_0,V_{1})$.
\item[\rm(2)] let $\sigma:\g_0\wedge V_0\longrightarrow \g_1$ be a linear map satisfying $\sigma(\varphi_0 u,v)=\sigma(u,\varphi_0 v), \forall u,v\in V_0$. Then $\omega^\sigma$ defined by
\begin{eqnarray*}
\omega^\sigma(X,Y)(x,y)\triangleq X\sigma(Yx,y)+X\sigma(x,Yy)-Y\sigma(Xx,y)-Y\sigma(x,Xy)
\end{eqnarray*}
is a $2$-cocycle of $\La(\g,\mathbb{V})$ with values in $\Hom(\wedge^2\frkg_0,V_{1})$.
\item[\rm(3)] 
    $C^1(\g,\mathbb{V})$ is a Lie algebra with the bracket $\{\cdot,\cdot\}$ given by
\[\{X+\xi,Y+\eta\}=[X,Y]_\varphi+X\triangleright \eta-Y\triangleright \xi+\omega^\sigma(X,Y).\]
\end{itemize}
\end{lem}
\pf To prove that $\triangleright$ is an action, we need to check
$$[Y,X]_\varphi \triangleright \xi=Y\triangleright(X\triangleright \xi)-X\triangleright(Y\triangleright \xi).$$
It follows from
\begin{eqnarray*}
&&\big(Y\triangleright(X\triangleright \xi)-X\triangleright(Y\triangleright \xi)\big)(x,y)
\\&=& Y\varphi (X\triangleright \xi)(x,y)-X\triangleright \xi(\varphi Y x,y)-X\triangleright \xi(x,\varphi Y y)-c.p.(X,Y)\\&=&
Y\varphi(X\varphi \xi(x,y)-\xi(\varphi Xx,y)-\xi(x,\varphi Xy))-X\varphi \xi(\varphi Yx,y)+\xi(\varphi X\varphi Yx,y)+\xi(\varphi Yx,\varphi X y)\\ &&-X\varphi \xi(x,\varphi Yy)+\xi(\varphi Xx,\varphi Yy)+\xi(x,\varphi X \varphi Yy)-c.p.(X,Y)\\ &=&[Y,X]_\varphi \varphi \xi(x,y)-\xi(\varphi[Y,X]_\varphi x,y)-\xi( x,\varphi[Y,X]_\varphi y)\\ &=&([Y,X]_\varphi \triangleright \xi)(x,y).
\end{eqnarray*}
It is direct to verify that $\omega^\sigma$ is a $2$-cocycle. We omit the details.\qed\vspace{3mm}
\subsection{Derivations of crossed modules}
The main goal of this section is to propose a notion of derivations of a crossed module of Lie $2$-algebras and prove that it is a strict Lie $2$-algebra, which generalises the derivations of a Lie $2$-algebra $\Der(\g)$.

Given a crossed module
$(\m,\g,\phi,\varphi,\sigma)$, since $\m$ is a $\g$-module, we have a natural
complex
\begin{equation*}
\CD
   \Der(\g,\m) : \Der_{1}(\g,\m)\triangleq \Hom(\g_0,\m_1)@>-\D>>\Der_{0}(\g,\m),
\endCD
\end{equation*}
where $\Der_{0}(\g,\m)$ is the set of $1$-cocycles and $\D$ is the Lie $2$-algebra coboundary operator. Denote by $\Inn_0(\g,\m)$ the set of $1$-coboundaries.
Explicitly, by (\ref{1-cob1}) and $(\ref{1-cob2}),$ we have $-\D(\Theta)=\delta(\Theta)+l_{\delta(\Theta)}$ and
\[l_{\delta(\Theta)}(x,y)={\Theta}[x,y]-x\triangleright \Theta{y}+y\triangleright {\Theta}x,\ \ \ \ \ \forall x,y\in \g_0,\Theta\in{\Der_{1}(\g,\m)}.\] We call $\Der(\g,\m)$ the {\bf derivations of the crossed module}.

By condition $(iv)$ of Definition \ref{defi:Lie 2 cm} and Lemma \ref{lem:degree 1 Lie alg}, $(C^1(\g,\m),\{\cdot,\cdot\})$ is a Lie algebra. As its subspace, $\Der_{0}(\g,\m)$ inherits a bracket operation, that is,
for any $X+l_{X}, Y+l_{Y}\in{\Der_{0}(\g,\m)}, $
\[\{X+l_{X},Y+l_{Y}\}=[X,Y]_\varphi+X\triangleright l_Y-Y\triangleright l_X+\omega^\sigma(X,Y).\]
Set $l_{[X,Y]_\varphi}=X\triangleright
l_Y-Y\triangleright l_X+\omega^\sigma(X,Y)$. To be more precise,
$$l_{[X,Y]_\varphi}(x,y)=X\varphi l_Y(x,y)-l_Y(\varphi Xx,y)-l_Y(x,\varphi Xy)+X\sigma (Yx,y)+X\sigma (x,Yy)-c.p.(X,Y).$$
Define $\{X+l_{X},\Theta\}\triangleq[X,\Theta]_\varphi.$

\begin{lem}\label{lem:sub Lie algebra}
With the above notations, $\Der_0(\g,\m)$ is a Lie subalgebra of
$(C^1(\g,\m),\{\cdot,\cdot\})$.
\end{lem}
\pf For any $X+l_{X}, Y+l_{Y}\in{\Der_{0}(\g,\m)},$ we need to prove that $[X,Y]_\varphi+l_{[X,Y]_\varphi}$ is a
$1$-cocycle. It is not hard to check conditions (\ref{d})-(\ref{b}).
We just verify condition (\ref{c}), i.e., for any
$x,y,z\in\g_0$,
\begin{eqnarray}\label{eq:d}
l_{[X,Y]_\varphi}(x,[y,z])+x\triangleright
l_{[X,Y]_\varphi}(y,z)-(y,z)\triangleright [X,Y]_\varphi x
+c.p.(x,y,z)=[X,Y]_\varphi l_3(x,y,z).
\end{eqnarray}
Firstly, as $X+l_X,Y+l_Y\in \Der_0(\g,\m)$ and
$\Pi=Id+\sigma+\varphi$ is a homomorphism, we have
\begin{eqnarray}\label{eq:b}
\nonumber\varphi X[y,z]&=&\varphi(y\triangleright Xz-z\triangleright
Xy+\tilde{d}l_X(y,z))\\ &=& [y,\varphi Xz]+d\sigma(y,Xz)-[z,\varphi
Xy]-d\sigma(z,Xy)+\varphi \tilde{d}l_X(y,z),
\end{eqnarray}
\begin{eqnarray}\label{eq:a}
\nonumber
X\varphi(Yl_3(x,y,z))&=&X\varphi\big(l_Y(x,[y,z])+x\triangleright
l_Y(y,z)-(y,z)\triangleright Yx +c.p.(x,y,z)\big)\\
\nonumber&=&X\varphi l_Y(x,[y,z])+x\triangleright X\varphi
l_Y(y,z)-\varphi l_Y(y,z)\triangleright Xx+l_X(x,d\varphi
l_Y(y,z))\\
&&+X\sigma(x,\tilde{d}l_Y(y,z))-X\varphi((y,z)\triangleright
Yx)+c.p.(x,y,z),
\end{eqnarray}
and
\begin{eqnarray}\label{eq:c}
\nonumber&&-X\big(\varphi((y,z)\triangleright Yx)+l_3(\varphi
Yx,y,z)\big)\\ \nonumber&=&
X\big(\sigma(Yx,[y,z])+\sigma(y,z\triangleright
Yx)-\sigma(z,y\triangleright Yx)
+[y,\sigma(z,Yx)]+[z,\sigma(Yx,y)]\big)\\
\nonumber&=&X\sigma(Yx,[y,z])+X\sigma(y,z\triangleright
Yx)-X\sigma(z,y\triangleright Yx)+y\triangleright
X\sigma(z,Yx)-\sigma(z,Yx)\triangleright Xy\\
&&+l_X(y,d\sigma(z,Yx))-z\triangleright
X\sigma(y,Yx)+\sigma(y,Yx)\triangleright Xz-l_X(z,d\sigma(y,Yx)).
\end{eqnarray}
Thus, $(\ref{eq:d})$ follows from
\begin{eqnarray*}
\nonumber&&l_{[X,Y]_\varphi}(x,[y,z])+x\triangleright
l_{[X,Y]_\varphi}(y,z)-(y,z)\triangleright [X,Y]_\varphi x
+c.p.(x,y,z)\\ \nonumber&=&X\varphi l_Y(x,[y,z])-l_Y(x,[y,\varphi
Xz]+d\sigma(y,Xz)-[z,\varphi Xy]-d\sigma(z,Xy)+\varphi
\tilde{d}l_X(y,z))\\ \nonumber &&-l_Y(\varphi X
x,[y,z])+X\sigma(Yx,[y,z])+X\sigma(x,y\triangleright
Yz-z\triangleright Yy+\tilde{d}l_Y(y,z))\\ \nonumber
&&+x\triangleright \big(X\varphi l_Y(y,z)-l_Y(\varphi
Xy,z)-l_Y(y,\varphi Xz)+X\sigma(Yy,z)+X\sigma(y,Yz)\big)\\ \nonumber
&&-(y,z)\triangleright X\varphi Yx-c.p.(X,Y)+c.p.(x,y,z)\ \ \ \ \
by\ (\ref{eq:b})\\ \nonumber &=& \big(X\varphi
l_Y(x,[y,z])+x\triangleright X\varphi l_Y(y,z)+l_X(x,d\varphi
l_Y(y,z)) +X\sigma(x,\tilde{d}l_Y(y,z))\big)\\ \nonumber&&
-\big(l_Y(\varphi Xx,[y,z])+c.p.(\varphi Xx,y,z)+y\triangleright
l_Y(z,\varphi Xx)+z\triangleright l_Y(\varphi
Xx,y)-(y,z)\triangleright Y\varphi Xx\big)\\
\nonumber&&+\big(X\sigma(Yx,[y,z])+c.p.(Yx,y,z) +y\triangleright
X\sigma(z,Yx)+l_X(y,d\sigma(z,Yx))
\\ \nonumber&&-z\triangleright X\sigma(y,Yx)-l_X(z,d\sigma(y,Yx))\big)-c.p.(X,Y)+c.p.(x,y,z)
\\ \nonumber&=&X\varphi Y l_3(x,y,z)+\big(\varphi l_Y(y,z)\triangleright Xx+X\varphi((y,z)\triangleright Yx)\\ \nonumber&&-Yl_3(\varphi Xx,y,z)+\varphi Xx\triangleright l_Y(y,z)-(z,\varphi Xx)\triangleright Yy-(\varphi Xx,y)\triangleright Yz\\ \nonumber&&-X\varphi ((y,z)\triangleright Yx)-Xl_3(\varphi Yx,y,z)+\sigma(z,Yx)\triangleright Xy-\sigma(y,Yx)\triangleright Xz+c.p.(x,y,z)\big)\\ \nonumber&&-c.p.(X,Y)\ \ \ \ \ by\ (\ref{eq:a}),\ (\ref{c})\ and\ (\ref{eq:c})\\ \nonumber&=&
l_{\phi_0(z)}(Xx,Yy)-l_{\phi_0(y)}(Xx,Yz)+c.p.(x,y,z)-c.p.(X,Y)+[X,Y]_\varphi
l_3(x,y,z)\\ &=&[X,Y]_\varphi l_3(x,y,z),
\end{eqnarray*}
where the penultimate equation holds since conditions $(i)$ and
$(iii)$ of Definition \ref{defi:Lie 2 cm}. \qed\vspace{3mm}
\begin{thm}\label{thm:Der(g,m)}
For a crossed module $(\m,\g,\phi,\varphi,\sigma)$, its derivations
$(\Der(\g,\m),\{\cdot,\cdot\})$ is a strict Lie $2$-algebra.
\end{thm}
\pf By Lemma \ref{lem:END} and Lemma \ref{lem:sub Lie algebra}, we
only need to verify that $-\D$ is a graded derivation with respect
to the bracket operation $\{\cdot,\cdot\}$, i.e., for any
$X+l_X\in{\Der_{0}(\g,\m)}, \Theta, \Theta' \in{\Der_{1}(\g,\m)}$,
\begin{eqnarray}
\label{eq:bracket9}-\D\{X+l_X,\Theta\}&=&\{X+l_X,-\D\Theta\},\\
\label{eq:bracket10} \{-\D\Theta,\Theta'\}&=&\{\Theta,-\D\Theta'\}.
\end{eqnarray}
The left hand side of (\ref{eq:bracket9}) is equal to
$\delta([X,\Theta]_\varphi)+l_{\delta([X,\Theta]_\varphi)}.$ By repeatedly applying condition $(i)$ of Definition \ref{defi:Lie 2
cm} and the fact that $\Pi$ is a homomorphism, we have
\begin{eqnarray*}
&&l_{\delta([X,\Theta]_\varphi)}(x,y)\\
&=&[X,\Theta]_\varphi[x,y]-x\triangleright [X,\Theta]_\varphi y
+y\triangleright [X,\Theta]_\varphi x\\&=&X\varphi (x\triangleright
\Theta{y}-y\triangleright {\Theta}x+l_{\delta(\Theta)}(x,y))-\Theta
\varphi(x\triangleright Xy-y\triangleright Xx+\tilde{d}l_X(x,y))\\
&&-x\triangleright(X\varphi \Theta y-\Theta\varphi X
y)+y\triangleright(X\varphi \Theta x-\Theta\varphi X
x)\\&=&X\big([x,\varphi\Theta{y}]+\sigma(x,\tilde{d}\Theta{y})-
[y,\varphi\Theta{x}]-\sigma(y,\tilde{d}\Theta{x}))+\varphi
l_{\delta(\Theta)}(x,y)\big)\\&& -\Theta\big([x,\varphi
X{y}]+d\sigma(x,Xy)-[y,\varphi X{x}]-d\sigma(y,Xx)+\varphi
\tilde{d}l_X(x,y)\big)\\ &&-x\triangleright(X\varphi \Theta
y-\Theta\varphi X y)+y\triangleright(X\varphi \Theta x-\Theta\varphi
X x)\\&=& x\triangleright X\varphi\Theta y-\varphi\Theta
y\triangleright Xx+l_X(x,d\varphi\Theta
y)+X\sigma(x,\tilde{d}\Theta{y})
\\&&-y\triangleright X\varphi\Theta x+\varphi\Theta x\triangleright Xy-l_X(y,d\varphi\Theta x)-X\sigma(y,\tilde{d}\Theta{x})+X\varphi l_{\delta(\Theta)}(x,y)
\\ &&-x\triangleright \Theta \varphi Xy+\varphi Xy\triangleright \Theta x-l_{\delta(\Theta)}(x,\varphi Xy)-\Theta d\sigma(x,Xy)
\\ &&+y\triangleright \Theta \varphi Xx-\varphi Xx\triangleright \Theta y+l_{\delta(\Theta)}(y,\varphi Xx)+\Theta d\sigma(y,Xx)-\Theta\varphi \tilde{d}l_X(x,y)
\\ &&-x\triangleright(X\varphi \Theta y-\Theta\varphi X y)+y\triangleright(X\varphi \Theta x-\Theta\varphi X x)
\\&=&l_X(x,d\varphi\Theta y)+X\sigma(x,\tilde{d}\Theta{y})-l_X(y,d\varphi\Theta x)-X\sigma(y,\tilde{d}\Theta{x})+X\varphi l_{\delta(\Theta)}(x,y)\\ &&-l_{\delta(\Theta)}(x,\varphi Xy)-\Theta d\sigma(x,Xy)+l_{\delta(\Theta)}(y,\varphi Xx)+\Theta d\sigma(y,Xx)-\Theta\varphi \tilde{d}l_X(x,y)\\ &=&l_{[X,\delta(\Theta)]_\varphi}(x,y).
\end{eqnarray*}
The right hand side
of (\ref{eq:bracket9}) is equal to
$$\{X+l_X,\delta(\Theta)+l_{\delta(\Theta)}\}=[X,\delta(\Theta)]_\varphi+
l_{[X,\delta(\Theta)]_\varphi}.$$ Thus, (\ref{eq:bracket9}) holds
since $\delta([X,\Theta]_\varphi)=[X,\delta(\Theta)]_\varphi.$

The equation (\ref{eq:bracket10}) is a consequence of
$[\delta(\Theta),\Theta']_\varphi=[\Theta,\delta(\Theta')]_\varphi$.
This finishes the proof.\qed\vspace{3mm}
\begin{pro}\label{pro:1st}
For a crossed module $(\m,\g,\phi,\varphi,\sigma)$, the 1st
cohomology group $\Ha^1(\g,\m)=\Der_0(\g,\m)/\Inn_0(\g,\m)$ is a
quotient Lie algebra.
\end{pro}
\pf We shall prove $\Inn_0(\g,\m)$ is an ideal of $\Der_0(\g,\m)$.
Since (\ref{eq:bracket9}) holds, it suffices to show
\begin{eqnarray}\label{not ideal}
\label{eq:bracket1}\{X+l_X,-\D\alpha\}&=&-\D(X\varphi\alpha+l_X(\varphi
\alpha,\cdot)+X\sigma(\alpha,\cdot)),\ \ \ \ \forall \alpha\in \m_0.
\end{eqnarray}
Acting on $x\in\g_0$ and taking into account $(\ref{1-cob1})$ and $(\ref{1-cob2})$, we
have
\begin{eqnarray*}
\{X+l_X,-\D\alpha\}(x)&=&-X\varphi(\D\alpha)x+(\D\alpha) \varphi
Xx\\ &=& -X\varphi(x\triangleright \alpha)+\varphi Xx\triangleright
\alpha\\ &=& -x\triangleright X\varphi \alpha+\varphi
\alpha\triangleright Xx-\tilde{d}l_X(x,\varphi
\alpha)-Xd\sigma(x,\alpha)+\varphi Xx\triangleright \alpha\\
&=&-x\triangleright X\varphi \alpha-\tilde{d}l_X(x,\varphi
\alpha)-\tilde{d}X\sigma(x,\alpha)\\
&=&-\D(X\varphi\alpha+l_X(\varphi
\alpha,\cdot)+X\sigma(\alpha,\cdot))(x),
\end{eqnarray*}
where we have used condition $(i)$ and $\Pi$ is a homomorphism in
Definition \ref{defi:Lie 2 cm}. Likewise, (\ref{eq:bracket1}) holds
on $\g_1$. Finally, relying on the coherence law of the
homomorphism $\Pi$, we have
\begin{eqnarray}\label{coherence}
\nonumber&&X\big(\sigma([x,y],\alpha)+ \sigma(y\triangleright
\alpha,x)-\sigma(x\triangleright
\alpha,y)-\varphi((x,y)\triangleright\alpha)-l_3(x,y,\varphi\alpha)\big)\\
\nonumber&=&X\big([x,\sigma(y,\alpha)]+[y,\sigma(\alpha,x)]\big)\\
\nonumber&=&x\triangleright
X\sigma(y,\alpha)-\sigma(y,\alpha)\triangleright
Xx+l_X(x,d\sigma(y,\alpha))\\ &&+y\triangleright
X\sigma(\alpha,x)-\sigma(\alpha,x)\triangleright
Xy+l_X(y,d\sigma(\alpha,x)).
\end{eqnarray}
Thus, acting on $\wedge^2 \g_0$, it turns out that
\begin{eqnarray*}
&&\{X+l_X,-\D\alpha\}(x,y)\\ &=&-X \varphi
\D\alpha(x,y)+\D\alpha(\varphi Xx,y)+\D\alpha(x,\varphi Xy)
+\D\alpha(\varphi l_X(x,y))-l_X(\varphi \D\alpha (x),y)\\
&&-l_X(x,\varphi \D\alpha
(y))-X\sigma(\D\alpha(x),y)-X\sigma(x,\D\alpha(y))+\D\alpha(\sigma(Xx,y))+\D\alpha(\sigma(x,Xy))\\
&=&-X \varphi((x,y)\triangleright \alpha)+(\varphi
Xx,y)\triangleright \alpha+(x,\varphi Xy)\triangleright
\alpha+\varphi l_X(x,y)\triangleright
\alpha-l_X([x,\varphi\alpha]+d\sigma(x,\alpha),y)\\
&&-l_X(x,[y,\varphi\alpha]+d\sigma(y,\alpha))-X\sigma(x\triangleright
\alpha,y)-X\sigma(x,y\triangleright \alpha)
+\sigma(Xx,y)\triangleright \alpha+\sigma(x,Xy)\triangleright
\alpha\\ &=& \big(x\triangleright
X\sigma(y,\alpha)-\sigma(y,\alpha)\triangleright Xx+y\triangleright
X\sigma(\alpha,x)-\sigma(\alpha,x)\triangleright
Xy-X\sigma([x,y],\alpha)+Xl_3(x,y,\varphi\alpha)\big)\\ &&
+\big(l_X(\varphi\alpha,[x,y])-(x,y)\triangleright
X\varphi\alpha-(y,\varphi\alpha)\triangleright
Xx-(\varphi\alpha,x)\triangleright Xy+x\triangleright
l_X(y,\varphi\alpha)\\ &&+y\triangleright l_X(\varphi\alpha,x)
-Xl_3(x,y,\varphi\alpha)\big)+l_{\phi_0(y)}(Xx,\alpha)-l_{\phi_0(x)}(Xy,\alpha)\
\ \ \ \ by\ (\ref{coherence})\ and\ (\ref{c})
\\ &=&x\triangleright X\sigma(y,\alpha)+y\triangleright X\sigma(\alpha,x)-X\sigma([x,y],\alpha)+l_X(\varphi\alpha,[x,y])
+x\triangleright l_X(y,\varphi\alpha)+y\triangleright
l_X(\varphi\alpha,x)
\\ &&-(x,y)\triangleright X\varphi\alpha+\big(l_{\phi_0(y)}(\alpha,Xx)-l_{\phi_0(x)}(\alpha,Xy)+l_{\phi_0(y)}(Xx,\alpha)-l_{\phi_0(x)}(Xy,\alpha)\big)\\ &=&-\D(X\varphi\alpha+l_X(\varphi \alpha,\cdot)+X\sigma(\alpha,\cdot))(x,y),
\end{eqnarray*} where we have used $(iii)$ of Definition \ref{defi:Lie 2 cm}. The proof is finished. \qed\vspace{3mm}

\begin{rmk}
 Let $\m=\g,i=Id$ in Example \ref{ideal Lie 2 cm}. Then Theorem \ref{thm:Der(g,m)} recovers the theorem that $\Der(\g)$ is a strict Lie $2$-algebra in \cite{Sheng1} and Proposition \ref{pro:1st} implies that $\Ha^1(\g)=\Der_0(\g)/\Inn_0(\g)$ is a quotient Lie algebra. Also, it justifies our definition of $\Inn_0(\g)$, while it is defined by $\Img (\D|_{\g_0})$ in \cite{Sheng1} which is not an ideal of $\Der_0(\g)$  by (\ref{not ideal}).
\end{rmk}
\begin{rmk}
It should be interesting to consider the strict $2$-group integrated
from the strict Lie $2$-algebra $\Der(\g,\m)$. In particular, for
$\Der(\g)$, we can explicate the strict $2$-group $\Aut(\g)$
consisted of all the automorphisms of Lie $2$-algebra $\g$.
Furthermore, for a strict Lie $2$-algebra $\g$ with corresponding
strict Lie $2$-groups $\mathcal{G}$, it is interesting to establish
the connection between $\Der(\g)$ and the automorphism $2$-group
$\Aut(\mathcal{G})$ (\cite{Norri}).
\end{rmk}
\section{Classification of strong crossed modules via $\Ha^3$}\label{section6}
\subsection{Crossed modules from extensions of short exact sequences}
Suppose $(\g,d,[\cdot,\cdot],l_3)$ is a Lie $2$-algebra and
$\mathfrak{k}\subset\g$ is an ideal. This gives rise to a short
exact sequence of Lie $2$-algebras
$$ 0\rightarrow\mathfrak{k}\hookrightarrow \g
\stackrel{\pi}\longrightarrow\mathfrak{h}\rightarrow0,$$ where
$\mathfrak{h}=\g/\mathfrak{k}$ and $\pi$ is the canonical
projection, which is a strong Lie $2$-algebra homomorphism.
Furthermore, given an $\mathfrak{h}$-module $\mathbb{V}$, it is natural that
$\mathbb{V}$ endows with a $\g$-module structure and then a trivial
$\mathfrak{h}$-module structure as follows: $$x\triangleright
u=\tilde{x} \triangleright u,\ \ (x,y)\triangleright
u=(\tilde{x},\tilde{y}) \triangleright u, \ \ \ \forall x,y\in
\g,u\in \mathbb{V},\tilde{x}=\pi(x),\tilde{y}=\pi(y).$$

The natural projection $\pi:\g\longrightarrow\mathfrak{h}$ induces a
map of cochain complexes:
\begin{equation}\label{eq:ext2}
\CD
   \cdots @> \D^\mathfrak{h} >> C^2(\mathfrak{h},\mathbb{V}) @>\D^\mathfrak{h} >> C^{3}(\mathfrak{h},\mathbb{V}) @> \D^\mathfrak{h} >> \cdots \\
  @. @V \pi^* VV @V \pi^*VV  @.  \\
  \cdots @> \D^{\g} >> C^2(\g,\mathbb{V}) @>\D^{\g}>> C^{3}(\g,\mathbb{V}) @> \D^{\g} >> \cdots.
\endCD
\end{equation}
In this section, we always suppose $\alpha,\beta,\gamma\in\mathfrak{k}_0, \xi\in\mathfrak{k}_1$, $u,v,w\in V_0, m\in V_1,$ and $x,y\in\g_0$, $a\in\g_1$.
\begin{lem}\label{lem:der}
For a $2$-cochain  $\lambda\in C^2(\g,\mathbb{V})$, the
following three statements are equivalent:
\begin{itemize}
\item[\rm(1)] $i_{e}(\D^\g \lambda)=0, \forall e\in \mathfrak{k}$;
\item[\rm(2)] 
 $\D^\g \lambda=\pi^*\theta$ for a $3$-cocycle $\theta\in C^3(\mathfrak{h},\mathbb{V})$;
\item[\rm(3)] The map $\phi^\lambda=(\phi^\lambda_0,\phi^\lambda_1,\phi^\lambda_2): \g\rightarrow \End(\mathfrak{k}\oplus_\lambda \mathbb{V})$ given below defines an action of $\g$ on $\mathfrak{k}\oplus_\lambda\mathbb{V}$,
 where $$\mathfrak{k}\oplus_\lambda \mathbb{V}:\mathfrak{k}_1\oplus V_1\stackrel{d^\lambda}\longrightarrow \mathfrak{k}_0\oplus V_0$$ is a $2$-vector space
 with $d^\lambda(\xi+m)=d\xi+d^\mathbb{V}m+\lambda_0(\xi)$ and
\begin{equation}\label{action by der}
\left\{\begin{array}{rcl} x\triangleright_\lambda(\alpha+u)&=&[x,\alpha]+\lambda_1(x,\alpha)+\tilde{x}\triangleright u,\\
 x\triangleright_\lambda(\xi+m)&=&[x,\xi]+\lambda_2(x,\xi)+\tilde{x}\triangleright m,\\
a\triangleright_\lambda(\alpha+u)&=&[a,\alpha]+\lambda_2(a,\alpha)+\tilde{a}\triangleright u,\\
(x,y)\triangleright_\lambda(\alpha+u)&=&-l_3(x,y,\alpha)-\lambda_3(x,y,\alpha)+(\tilde{x},\tilde{y})\triangleright
u.\end{array}\right.
\end{equation}
That is, $\phi^\lambda$ is a Lie $2$-algebra homomorphism.
\end{itemize}
\end{lem}
\pf It is obvious that $(1)\Leftrightarrow (2)$, so it suffices to
prove $(1)\Leftrightarrow (3)$. Suppose that $\lambda$ satisfies
$(1)$, we shall show that $\phi^\lambda_0(x)\in \End_0(\mathfrak{k}\oplus_\lambda \mathbb{V})$ and $\phi^\lambda$ is a Lie $2$-algebra homomorphism. Referring to (\ref{3-ex cochain}), we use subscripts
to distinguish the five components in $C^3(\g,\mathbb{V})$.

Firstly, the equality $\phi^\lambda_0(x)\circ
d^\lambda=d^\lambda\circ \phi^\lambda_0(x)$ holds since
\begin{eqnarray*}
&&(\phi^\lambda_0(x)\circ d^\lambda-d^\lambda\circ
\phi^\lambda_0(x))(\xi+m)\\ &=&
x\triangleright_\lambda(d\xi+d^\mathbb{V}m+\lambda_0(\xi))-d^\lambda([x,\xi]+\lambda_2(x,\xi)+\tilde{x}\triangleright
m)\\ &=&
[x,d\xi]+\lambda_1(x,d\xi)+\tilde{x}\triangleright(d^\mathbb{V}m+\lambda_0(\xi))-d[x,\xi]
-d^\mathbb{V}\lambda_2(x,\xi)-d^\mathbb{V}(\tilde{x}\triangleright
m)-\lambda_0[x,\xi]\\ &=& \lambda_1(x,d\xi)+\tilde{x}\triangleright
\lambda_0(\xi)-d^\mathbb{V}\lambda_2(x,\xi)-\lambda_0[x,\xi]\\ &=&
(\D^\g\lambda)_0(x,\xi)\\ &=&0,
\end{eqnarray*}
where we have used the fact that $\triangleright$ is an action. Analogously, we get
\begin{eqnarray*}
(\phi^\lambda_0(da)-\delta \phi^\lambda_1(a))(\alpha+u)&=&-(\D^\g\lambda)_0(\alpha,a)=0,\\
(\phi^\lambda_0(da)-\delta \phi^\lambda_1(a))(\xi+m)&=&-(\D^\g\lambda)_1(a,\xi)=0.
\end{eqnarray*}

Secondly, by an elementary computation, we obtain
\begin{eqnarray*} &&\big(\phi^\lambda_0[x,y]-[\phi^\lambda_0(x),\phi^\lambda_0(y)]-\delta\phi^\lambda_2(x,y)\big)(\alpha+u)\\ &=&
[[x,y],\alpha]+\lambda_1([x,y],\alpha)+\widetilde{[x,y]}\triangleright u-[x,[y,\alpha]]-\lambda_1(x,[y,\alpha])-\tilde{x}\triangleright(\lambda_1(y,\alpha)+\tilde{y}\triangleright u)\\ &&
+[y,[x,\alpha]]+\lambda_1(y,[x,\alpha])+\tilde{y}\triangleright(\lambda_1(x,\alpha)+\tilde{x}\triangleright u)
+dl_3(x,y,\alpha)\\ &&+d^\mathbb{V}\lambda_3(x,y,\alpha)-d^\mathbb{V}((\tilde{x},\tilde{y})\triangleright
u)+\lambda_0l_3(x,y,\alpha)\\ &=&\lambda_1([x,y],\alpha)+c.p.-\tilde{x}\triangleright\lambda_1(y,\alpha)+\tilde{y}\triangleright\lambda_1(x,\alpha)
+d^\mathbb{V}\lambda_3(x,y,\alpha)+\lambda_0l_3(x,y,\alpha)\\ &=&-(\D^\g\lambda)_2(x,y,\alpha)\\ &=&0,
\end{eqnarray*}
where we have used the general Jacobi identity of $\g$ and the fact that $\triangleright$ is an action. Likewise, we can deduce that
\begin{eqnarray*}
\big(\phi^\lambda_0[x,y]-[\phi^\lambda_0(x),\phi^\lambda_0(y)]-\delta\phi^\lambda_2(x,y)\big)(\xi+m)&=&
-(\D^\g\lambda)_3(x,y,\xi)=0,\\ \big(\phi^\lambda_1[x,a]-[\phi^\lambda_0(x),\phi^\lambda_1(a)]-\delta\phi^\lambda_2(x,da)\big)(\alpha+u)&=&
-(\D^\g\lambda)_3(x,a,\alpha)=0.
\end{eqnarray*}

Finally, for the coherence condition, following from the fact that
$\triangleright$ is an action and the coherence law of
$[\cdot,\cdot]$ and $l_3$, we have
\begin{eqnarray*}
&&\big([\phi^\lambda_0(x),\phi^\lambda_2(y,z)]-\phi^\lambda_2([x,y],z)+c.p.-\phi^\lambda_1l_3(x,y,z)\big)(\alpha+u)
\\ &=&-[x,l_3(y,z,\alpha)]-\lambda_2(x,l_3(y,z,\alpha))-\tilde{x}\triangleright \lambda_3(y,z,\alpha)+\tilde{x}\triangleright((\tilde{y},\tilde{z})\triangleright u)\\ &&+l_3(y,z,[x,\alpha])+\lambda_3(y,z,[x,\alpha])-(\tilde{y},\tilde{z})\triangleright \lambda_1(x,\alpha)-(\tilde{y},\tilde{z})\triangleright (\tilde{x}\triangleright u)\\ &&+l_3([x,y],z,\alpha)+\lambda_3([x,y],z,\alpha)-(\widetilde{[x,y]},\tilde{z})\triangleright u+c.p.(x,y,z)\\ &&-[l_3(x,y,z),\alpha]-\lambda_2(l_3(x,y,z),\alpha)- \widetilde{l_3(x,y,z)}\triangleright u\\ &=&-(\D^\g\lambda)_4(x,y,z,\alpha)\\ &=&0.
\end{eqnarray*}Thus $\phi^\lambda$ is a Lie $2$-algebra homomorphism.
From the process above, the other hand is obvious.\qed\vspace{3mm}

Define $\varphi:\mathfrak{k}\oplus_\lambda
\mathbb{V}\longrightarrow\g$ by $\varphi=i\oplus 0$, which is a chain map.
  Assume that $\lambda$ satisfies the condition $i_{e}(\D^\g \lambda)=0, \forall e\in \mathfrak{k}$. Then by a simple check, we find that $(\mathfrak{k}\oplus_\lambda
\mathbb{V},\g,\phi^\lambda,\varphi)$ satisfies all the conditions in Proposition \ref{pro:Lie 2 cm}. Thus, we obtain:
\begin{pro}\label{pro:surjection}
Suppose that $i_{e}(\D^\g \lambda)=0, \forall e\in \mathfrak{k}$.
Then
 $\varepsilon_\lambda=(\mathfrak{k}\oplus_\lambda \mathbb{V},\g,\phi^\lambda,\varphi)$ is a strong crossed module, where
 the Lie $2$-algebra structure on $\mathfrak{k}\oplus_\lambda \mathbb{V}$ is given by
 \begin{equation}\label{Lie 2 st}
\left\{\begin{array}{rcll}
{[\alpha+u,\beta+v]}_\lambda&=&[\alpha,\beta]+\lambda_1(\alpha,\beta),\\
{[\alpha+u,\xi+m]}_\lambda&=&[\alpha,\xi]+\lambda_2(\alpha,\xi),\\
l^\lambda_3(\alpha+u,\beta+v,\gamma+w)&=&l_3(\alpha,\beta,\gamma)+\lambda_3(\alpha,\beta,\gamma),\end{array}\right.
  \end{equation}
  and $l_{\phi^\lambda_0(x)}$ is defined by
 $$l_{\phi^\lambda_0(x)}(\alpha+u,\beta+v)=l_3(x,\alpha,\beta)+\lambda_3(x,\alpha,\beta).$$
\end{pro}
Next, we deal with the problem of deformation. As will seen in the next proposition, we get the same strong crossed module in the isomorphic sense if $\lambda$ is modified by a coboundary and an element in $\Img\pi^*$.

Suppose $\lambda\in C^2(\g,\mathbb{V})$ satisfying that $\D^\g \lambda=\pi^*\theta$ for a $3$-cocycle $\theta\in C^3(\mathfrak{h},\mathbb{V})$. Then for any $R\in C^2(\mathfrak{h},\mathbb{V})$ and $A\in C^1(\g,\mathbb{V})$,
we have $$\D^\g (\lambda+\D^\g A+\pi^* R)=\pi^*(\theta+\D^{\mathfrak{h}}R).$$ Note that $\varepsilon_{\lambda+\D^\g A+\pi^* R}=\varepsilon_{\lambda+\D^\g A}$ due to $\pi^* R|_{\mathfrak{k}}=0$. Then define
$F:\mathfrak{k}\oplus_{\lambda+\D^\g A}
\mathbb{V}\longrightarrow\mathfrak{k}\oplus_{\lambda} \mathbb{V}$
 by
\begin{equation*}
\left\{\begin{array}{rcl} F_0(\alpha+u)&=&\alpha+u+A_0 (\alpha),\\
F_1(\xi+m)&=&\xi+m+A_1(\xi),\\
F_2(\alpha+u,\beta+v)&=&A_2(\alpha,\beta),
\end{array}\right.
\end{equation*}
$G=Id:\g\longrightarrow\g$, and
$\tau:\g_0\wedge(\mathfrak{k}\oplus_{\lambda+\D^\g A}\mathbb{V})_0\longrightarrow(\mathfrak{k}\oplus_{\lambda}\mathbb{V})_1$
by $\tau(x,\alpha+u)=A_2(x,\alpha)$.
\begin{pro}\label{two}
The map $(F,G,\tau)$ is an isomorphism from $\varepsilon_{\lambda+\D^\g A+\pi^* R}$ to $\varepsilon_{\lambda}.$ We call it a {\bf gauge transformation}.
\end{pro}
\pf Firstly,
$F$ is a Lie $2$-algebra homomorphism due to the formulation of $\D^\g
A$. Indeed, the condition $F_0\circ d^{\lambda+\D^\g A}=d^{\lambda}\circ
F_1$ follows from
\begin{eqnarray*}
&&(F_0\circ d^{\lambda+\D^\g A}-d^{\lambda}\circ F_1)(\xi+m)\\ &=&F_0(d\xi+d^\mathbb{V}m+(\lambda_0+(\D^\g A)_0)(\xi))-d^{\lambda}(\xi+m+A_1(\xi))\\
&=&d\xi+d^\mathbb{V}m+\lambda_0(\xi)+(\D^\g A)_0(\xi)+A_0(d\xi)-d\xi-d^\mathbb{V}(m+A_1(\xi))-\lambda_0(\xi)
\\ &=&(\D^\g A)_0(\xi)+A_0(d\xi)-d^\mathbb{V}A_1(\xi)\\ &=&0.
\end{eqnarray*}
Then, we
shall verify the coherence condition, since the other two conditions of homomorphism
are similar to get. By a straightforward calculation, we have
\begin{eqnarray*}
&&[F_0(\alpha+u),F_2(\beta+v,\gamma+w)]_{\lambda}-F_2([\alpha+u,\beta+v]_{\lambda+\D^\g A},\gamma+w)+c.p.\\
&&+l^{\lambda}_3(F_0(\alpha+u),F_0(\beta+v),F_0(\gamma+w))-F_1l^{\lambda+\D^\g A}_3(\alpha+u,\beta+v,\gamma+w)
\\ &=&-A_2([\alpha,\beta],\gamma)+c.p.+l_3(\alpha,\beta,\gamma)
-l_3(\alpha,\beta,\gamma)-(\D^\g A)_3(\alpha,\beta,\gamma)-A_1l_3(\alpha,\beta,\gamma)\\ &=&0.
\end{eqnarray*} Next,  we get $\varphi\circ F=G\circ \varphi$ by definition.
Finally, similar to the above procedure, it is direct to check that
$$\big((G_0,F_0),(G_1,F_1),(G_2,\tau,F_2)\big):\g\triangleright_{\lambda+\D^\g A}(\mathfrak{k}\oplus_{\lambda+\D^\g A}\mathbb{V})\longrightarrow
\g\triangleright_{{\lambda}}(\mathfrak{k}\oplus_{\lambda}\mathbb{V})$$
is a Lie $2$-algebra homomorphism. Thus, $(F,G,\tau)$ is a morphism
from $\varepsilon_{\lambda+\D^\g A}$ to $\varepsilon_{\lambda}$. Furthermore,
note that $F$ and $G$ are bijections as chain maps and $\varepsilon_{\lambda+\D^\g A+\pi^* R}=\varepsilon_{\lambda+\D^\g A}$. We conclude that
$\varepsilon_{\lambda+\D^\g A+\pi^* R}$ and $\varepsilon_{\lambda}$ are
isomorphic.\qed\vspace{3mm}

Generally, for a $3$-cocycle $\theta\in
C^3(\mathfrak{h},\mathbb{V})$, the $3$-cochain $\pi^*\theta$ is not necessarily to
be a $3$-coboundary. Nevertheless, consider the short exact
sequence $$0\rightarrow\ker
\pi\hookrightarrow\mathfrak{F}(\mathfrak{h})\stackrel{\pi}\longrightarrow\mathfrak{h}\rightarrow0,$$
where  $\mathfrak{F}(\mathfrak{h})$ is the free Lie $2$-algebra
(\cite{Markl}) generated by the underlying $2$-vector space of
$\mathfrak{h}$ and $\pi$ is the canonical projection. For a
$3$-cocycle $\theta\in C^3(\mathfrak{h},\mathbb{V})$, since the
second and third cohomology of any free Lie $2$-algebra are
trivial\footnote{It is due to that the cohomology of Lie $2$-algebra
is ``operadic'' in nature and the operadic cohomology (degree $\geq
2$) vanishes on frees. See \cite{AQ}.}, there exists a $2$-cochain
$\lambda\in C^2(\mathfrak{F}(\mathfrak{h}),\mathbb{V})$ such that
$\D^{\mathfrak{F}(\mathfrak{h})}\lambda=\pi^*\theta$ and for
different $2$-cochains $\lambda,\lambda'$ satisfying it, we have
$[\lambda-\lambda']=0$. By Proposition \ref{two}, we
have:
\begin{cor}\label{main2}
For any $[\theta]\in \Ha^3(\mathfrak{h},\mathbb{V})$, we get a class
of crossed modules differing from each other by a gauge
transformation
$$\big\{\varepsilon_\lambda;\D^{\mathfrak{F}(\mathfrak{h})}\lambda=\pi^*\theta',
[\theta']=[\theta]\in\Ha^3(\mathfrak{h},\mathbb{V})\big\}.$$
\end{cor}
\subsection{Classification}
A crossed module $(\m,\g,\phi,\varphi,\sigma)$ can yield a 4-term
exact sequence of 2-vector spaces
\begin{equation}\label{seq:Lie 2 cm}
\xymatrix@C=0.5cm{ 0 \ar[r] & \mathbb{V }\ar[rr]^{i} &&
                \mathfrak{m} \ar[rr]^{\varphi} && \mathfrak{g} \ar[rr]^{\pi} && \mathfrak{h}\ar[r]  & 0,
                }
\end{equation}
where $\mathbb{V}\triangleq\ker\varphi, \mathfrak{h}\triangleq\coker
\varphi$, and $i,\pi$ are the canonical inclusion and projection. By
$(i),(ii)$ of Definition \ref{defi:Lie 2 cm}$, \mathbb{V}$ is in the
center of $\mathfrak{m}$.
However, it needs some extra conditions to
ensure that $\mathfrak{h}$ is a Lie $2$-algebra and there exists an
induced action of $\mathfrak{h}$ on $\mathbb{{V}}$.
\begin{lem}\label{lem:induced action}
For a crossed module of Lie $2$-algebras
$(\m,\g,\phi,\varphi,\sigma)$, we have
\begin{itemize}
\item[\rm(1)] $\mathfrak{h}=\mathfrak{g}/\Img\varphi$ is a quotient Lie $2$-algebra of $\mathfrak{g}$ if $\Img\sigma\subset{\Img\varphi_1}$;
\item[\rm(2)] the action of $\mathfrak{g}$ on $\mathfrak{m}$ induces an $\mathfrak{h}$-module structure on $\mathbb{V}$ if $\Img\sigma\subset{\Img\varphi_1}$ and $\sigma(\ker\varphi_{0}\wedge \g_0)=0.$
\end{itemize}
In particular, $\mathbb{V}$ becomes a module of Lie $2$-algebra
$\mathfrak{h}$ if the crossed module is strong.
\end{lem}
\pf
In order to prove $\Img\varphi$ is an ideal, we need to verify
$l_2(\Img\varphi\wedge \mathfrak{g})\subset{\Img\varphi},
l_3(\Img\varphi_0\wedge \mathfrak{g}_{0}\wedge
\mathfrak{g}_{0})\subset{\Img\varphi_1}.$ Following that
$\Pi=Id+\sigma+\varphi$ is a Lie 2-algebra homomorphism, we have
\begin{equation}\label{homo}
\left\{\begin{array}{rcll} \varphi_0(x\triangleright\alpha)-[x,\varphi_0(\alpha)]&=&d\sigma(x,\alpha),&~\forall x\in\g_0,\alpha\in\m_0,\\
\varphi_1(a\triangleright \alpha)-[a,\varphi_0(\alpha)]&=&\sigma(da,\alpha),&~\forall a\in\g_1,\alpha\in\m_0,\\
\varphi_1(x\triangleright\xi)-[x,\varphi_1(\xi)]&=&\sigma(x,\tilde{d}\xi),&~\forall
x\in\g_0,\xi\in\m_1,
\end{array}\right.
\end{equation}
which leads to $l_2(\Img\varphi\wedge
\mathfrak{g})\subset{\Img\varphi}$ since $\Img\sigma
\subset{\Img\varphi_{1}}$. Moreover, consider the coherence law in
the definition of homomorphism, for any
$x,y\in{\mathfrak{g}_0},\alpha\in{\mathfrak{m}_0}$,
\[[\sigma(y,\alpha),x]+[\sigma(\alpha,x),y]+\sigma([x,y],\alpha)+
\sigma(y\triangleright \alpha,x)-\sigma(x\triangleright
\alpha,y)=l_3(x,y,\varphi_0\alpha)-\varphi_1L_3(x,y,\alpha),\]
which implies that $l_3(x,y,\varphi_0\alpha)\in{\Img\varphi_1},$
i.e., $l_3(\Img\varphi_0\wedge \mathfrak{g}_{0}\wedge
\mathfrak{g}_{0})\subset{\Img\varphi_1}.$

To prove $(2)$, firstly, by equalities (\ref{homo}) and
$\sigma(\ker\varphi_{0}\wedge \g_0)=0$, we get
$\phi_i:\mathfrak{g_i}\longrightarrow{\End_i(\ker\varphi)},i=0,1.$
Coupled with the coherence law above, it is clear that
\[(x,y)\triangleright \beta=-L_3(x,y,\beta)\in \ker\varphi_1,\ \ \forall \beta \in \ker\varphi_0, x,y\in \g_0.\] Namely, $\phi_2:\wedge^{2}\mathfrak{g}_{0}\longrightarrow{\End_1(\ker\varphi)}.$ Thus we can define $\tilde{\phi}=\phi\circ s: \mathfrak{h}\longrightarrow\End(\ker\varphi)$, where $s:\mathfrak{h}\longrightarrow\g$ is a section of $\pi$.

Next, we prove $\tilde{\phi}$ is independent of the section $s$. Let $s'$ be another section of $\pi$, then $\Img(s-s')\in\ker \pi=\Img \varphi$. We shall check $\phi\circ (s-s')=0$, which is equivalent to show the homomorphism $\g\stackrel{\phi}\longrightarrow \End(\ker\varphi)$ vanishes when restricting on $\Img\varphi$. Depending on $(i),(iii)$ of Definition \ref{defi:Lie 2 cm}, we have
\[\varphi(\alpha)\triangleright{\beta}=-\varphi(\beta)\triangleright{\alpha}=0,\ \ \forall{\beta}\in{\ker\varphi},\alpha\in \m.\]
And, for any ${\beta}\in{\ker\varphi_0}, \alpha\in \m_0, x\in \g_0$,
\[(x,\varphi_{0} \alpha)\triangleright
\beta=-(x,\varphi_{0}\beta)\triangleright
\alpha-\sigma(x,\beta)\triangleright \alpha-
\sigma(x,\alpha)\triangleright \beta=-
\sigma(x,\alpha)\triangleright \beta,\] which vanishes since
$\Img\sigma \subset \Img\varphi_1$ and $\beta\in{\ker\varphi_0}$.

At last, since $[sx,sy]\triangleright\alpha=s[x,y]\triangleright\alpha, \forall x,y\in \mathfrak{h},\alpha \in \ker \varphi,$ it is direct to check $\tilde{\phi}$ is a Lie $2$-algebra homomorphism.
Therefore, the requirements $\Img\sigma \subset \Img\varphi_1$ and
$\sigma(\ker\varphi_{0}\wedge \g_0)=0$ make sure that the action of
$\mathfrak{g}$ on $\mathfrak{m}$ can induce an action of
$\mathfrak{h}$ on $\ker\varphi.$\qed\vspace{3mm}


Note that for a strong crossed module $(\m,\g,\phi,\varphi)$, there exists an exact sequence (\ref{seq:Lie 2 cm}) satisfying that all the homomorphisms are strong. Denote by $\mathcal{C}(\mathfrak{h},\mathbb{V})$ the set of strong crossed modules with respect to fixed kernel $\mathbb{V}$ and fixed cokernel $\mathfrak{h}$ and fixed action of $\mathfrak{h}$ on $\mathbb{V}$.
\begin{ex}{\rm
 Consider the set in Proposition \ref{main2} denoted by $\mathcal{C}(\mathfrak{h},\mathbb{V})_{[\theta]}$. In fact,
$$\mathcal{C}(\mathfrak{h},\mathbb{V})_{[\theta]}=\big\{\varepsilon_\lambda;\D^{\mathfrak{F}(\mathfrak{h})}\lambda=\pi^*\theta',
[\theta']=[\theta]\in\Ha^3(\mathfrak{h},\mathbb{V})\big\}$$ is a
subset of $\mathcal{C}(\mathfrak{h},\mathbb{V})$.}
\end{ex}

\begin{defi}\label{defi:el}
Let $(\m,\g,\phi,\varphi)$ and $(\m',\g',\phi',\varphi')$ be two
strong crossed modules in  $\mathcal{C}(\mathfrak{h},\mathbb{V})$, a
strong map between them is a strong morphism of crossed modules $(F,G)$
such that
the diagram
\begin{equation}\label{eq:ext1}
\CD
  0 @> >>  \mathbb{V} @>i>> \mathfrak{m} @>\varphi>> \mathfrak{g}@>\pi>> \mathfrak{h} @> >> 0 \\
 @. @V Id VV @V FVV @VG VV @V Id VV @.  \\
 0 @> >> \mathbb{V} @>i'>> \mathfrak{m'} @>\varphi'>> \mathfrak{g'}@>\pi'>> \mathfrak{h} @> >> 0,
\endCD
\end{equation}
is commutative.
\end{defi}
\begin{defi}\label{equi re}
For two crossed module $\varepsilon,\varepsilon'\in
\mathcal{C}(\mathfrak{h},\mathbb{V})$, we define $\varepsilon\sim\varepsilon'$
 if there exist two
crossed module $\varepsilon_\lambda, \varepsilon_{\lambda'}\in
\mathcal{C}(\mathfrak{h},\mathbb{V})_{[\theta]}$ for a $3$-cocycle $\theta\in C^3(\mathfrak{h},\mathbb{V})$ such that the
diagram
\begin{equation}\label{diag}
\xymatrix{
                & \varepsilon_\lambda\cong \varepsilon_{\lambda'} \ar[dl]_{f} \ar[dr]^{f'}    \\
   \varepsilon  &\sim &\varepsilon'   }
\end{equation}
holds, where $f, f'$ are strong maps in
$\mathcal{C}(\mathfrak{h},\mathbb{V})$ and $\cong$ is the gauge transformation between them.
\end{defi}
Since the composition of gauge transformations is also a gauge
transformation, it is evident that $\sim$ is an equivalence
relation.
\begin{rmk}
The linear map $\varphi:\m\rightarrow \g$ induces an action groupoid
$\g\times \m\rightrightarrows\g$ with the abelian group structure on
$\m$, where the source, target and inclusion maps are
$s(x,\alpha)=x, t(x,\alpha)=x+\varphi(\alpha), i(x)=(x,0).$ Then
$\varphi$ has fixed kernel and cokernel means that the groupoid has
fixed isotropy subgroup and orbit space. In this viewpoint, the
relation between $\varepsilon$ and $\varepsilon'$ defined by
(\ref{diag}) is in fact a generalized map (\cite{gm}) between them.
\end{rmk}
\begin{thm}\label{thm:main}For a Lie $2$-algebra $\mathfrak{h}$ and an $\mathfrak{h}$-module $\mathbb{V}$, there exists a canonical bijection$$\mathcal{C}(\mathfrak{h},\mathbb{V})/_\sim \stackrel{\approx}\longrightarrow\Ha^{3}(\mathfrak{h},\mathbb{V}),$$
where $\sim$ is the equivalence relation defined in Definition
\ref{equi re}.
\end{thm}

We divide the proof into four steps. For simplicity, we denote by
the same notations $d,[\cdot,\cdot]$ and $l_3$ for different Lie
$2$-algebras except when emphasis is needed, which will not cause
any confusion. In the rest of this section, we will always suppose
$x,y,z\in{\mathfrak{h}_{0}}$ and $a,b\in{\mathfrak{h}_{1}}.$

Step $1$: Construct a map $\mu: \mathcal{C}(\mathfrak{h},\mathbb{V})\longrightarrow\Ha^{3}(\mathfrak{h},\mathbb{V})$. 
To a crossed module $\varepsilon=(\m,\g,\phi,\varphi)\in \mathcal{C}(\mathfrak{h},\mathbb{V})$,
choose linear sections $s=(s_0,s_1):\mathfrak{h}\rightarrow \g,\ \pi
s=Id$ and $q=(q_0,q_1):\Img \varphi\rightarrow \m,\ \varphi q=Id$.
Since $ds_1(a)-s_{0}d(a)\in \ker \pi_0=\Img \varphi_0$, we can take
$$\lambda_0(a)= q_0(ds_1(a)-s_{0}d(a)).$$ Similarly, take
  \begin{equation*}
\left\{\begin{array}{rcl}
\lambda_1(x,y)&=& q_0([s_0(x),s_0(y)]-s_0[x,y]),\\
\lambda_2(x,a)&=& q_1([s_0(x),s_1(a)]-s_1[x,a]),\\
\lambda_3(x,y,z)&=&
q_1(l_3(s_0(x),s_0(y),s_0(z))-s_{1}l_3(x,y,z)).\end{array}\right.
  \end{equation*}
Note that $\lambda_\varepsilon=\sum_{i=0}^{3}\lambda_i$ satisfies
$\pi^*\lambda_\varepsilon\in C^2(\g,\m)$. It is reasonable to define
\begin{eqnarray}\label{34}
\theta_\varepsilon=s^*\D^\g(\pi^*\lambda_\varepsilon).
\end{eqnarray}

\begin{lem}\label{lem:3-cocycle}With the above notations, we have
\begin{itemize}
\item[\rm(1)] $\varphi\circ \theta_\varepsilon=0,$ that is, $\theta_\varepsilon\in C^3(\mathfrak{h},\mathbb{V}).$
\item[\rm(2)] $\D^\mathfrak{h}\theta_\varepsilon=0.$
\end{itemize}
\end{lem}
\pf According to $(\ref{3-ex cochain})$, we shall check that
$\varphi_0\circ \theta_{\varepsilon j}=0$ for j=0,2 and
$\varphi_1\circ \theta_{\varepsilon j}=0$ for j=1,3,4. The cases of
$j=0,1,2,3$ are quite straightforward, so we omit the details.
For the case of $j=4$, since $\Pi=Id+\varphi$ is a strong
homomorphism, we have
\begin{eqnarray}\label{eq:bracket15}
l_3(s_0(x_{\sigma_1}),s_0(x_{\sigma_2}),\varphi_0\lambda_1(x_{\sigma_3},x_{\sigma_4}))=
-\varphi_1((s_0(x_{\sigma_1}),s_0(x_{\sigma_2}))\triangleright
\lambda_1(x_{\sigma_3},x_{\sigma_4})).
\end{eqnarray}
 Taking into account $(\ref{3-ex cochain}),(\ref{eq:bracket15})$ and the coherence laws of $l_2,l_3$ of Lie 2-algebras $\mathfrak{g}$ and $\mathfrak{h}$, we have
\begin{eqnarray*}
&&\varphi_1(\theta_{\varepsilon 4}(x_1,x_2,x_3,x_4))\\
&=&-\sum_{\sigma}(-1)^{\sigma}l_3(s_0(x_{\sigma_1}),
s_0(x_{\sigma_2}),[s_0(x_{\sigma_3}),s_0(x_{\sigma_4})]-s_0[x_{\sigma_3},x_{\sigma_4}])\\
&&-\sum_{\tau}(-1)^{\tau}\{[s_0(x_{\tau_4}),s_{1}l_3(x_{\tau_1},x_{\tau_2},x_{\tau_3})]-s_1[x_{\tau_4},
l_3(x_{\tau_1},x_{\tau_2},x_{\tau_3})]\}
\\ &&+\sum_{i=1}^{4}(-1)^{i+1}[s_0(x_i),l_3(s_0(x_1),\cdots,s_0(\widehat{x_i}),\cdots,s_0(x_{4}))
-s_1l_3(x_1,\cdots,\widehat{x_i},\cdots,x_{4})] \\ &&
+\sum_{i<j}(-1)^{i+j}\{l_3(s_0[x_i,x_j],s_0(x_1),\cdots,s_0(\widehat{x_i}),\cdots,s_0(\widehat{x_j}),\cdots,s_0(x_{4}))\\
&&-
s_1l_3([x_i,x_j],x_1,\cdots,\widehat{x_i},\cdots,\widehat{x_j},\cdots,x_{4})\}\\
&=&-\sum_{\sigma}(-1)^{\sigma}l_3(s_0(x_{\sigma_1}),
s_0(x_{\sigma_2}),[s_0(x_{\sigma_3}),s_0(x_{\sigma_4})])
\\ &&+\sum_{i=1}^{4}(-1)^{i+1}[s_0(x_i),l_3(s_0(x_1),\cdots,s_0(\widehat{x_i}),\cdots,s_0(x_{4}))]\\ &&+\sum_{\tau}(-1)^{\tau}s_1[x_{\tau_4},
l_3(x_{\tau_1},x_{\tau_2},x_{\tau_3})] +\sum_{i<j}(-1)^{i+j}
s_1l_3([x_i,x_j],x_1,\cdots,\widehat{x_i},\cdots,\widehat{x_j},\cdots,x_{4})\\
&=&0.
\end{eqnarray*}
Therefore, $\Img \theta_\varepsilon\subset{\mathbb{V}}$.

Note that if only $\Img \lambda_\varepsilon\subset \mathbb{V}$, the
definition of $\theta_\varepsilon$ would be read as
$\theta_\varepsilon=s^*\pi^*\D^\mathfrak{h}(\lambda_\varepsilon)=\D^\mathfrak{h}\lambda_\varepsilon$
and hence would give $\D^\mathfrak{h}\theta_\varepsilon=0$. In our
situation, straightforward calculations show that
$\D^\mathfrak{h}\theta_\varepsilon$ still vanishes. Explicitly, by
the definition of $\D^\mathfrak{h}$, we have
 \begin{equation*}
\left\{\begin{array}{rcll}
(\D^\mathfrak{h}\theta_\varepsilon)_0&=&\hat{d}\theta_0+\widehat{d^\mathbb{V}}\theta_1&~\in\Hom(\odot^2\g_1,V_0),\\
(\D^\mathfrak{h}\theta_\varepsilon)_1&=& d^{(1,0)}_{\phi}\theta_0+\hat{d}\theta_2+\widehat{d^\mathbb{V}}\theta_3&~\in\Hom(\wedge^2\g_0\wedge\g_1,V_0),\\
(\D^\mathfrak{h}\theta_\varepsilon)_2&=&d^{(0,1)}_{\phi}\theta_0+d^{(1,0)}_{\phi}\theta_1+\hat{d}\theta_3&~\in\Hom(\g_0\wedge\odot^2\g_1,V_1),\\
(\D^\mathfrak{h}\theta_\varepsilon)_3&=&d_{\phi_2}\theta_0+d_{l_3}\theta_1+d^{(0,1)}_{\phi}\theta_2+d^{(1,0)}_{\phi}\theta_3+\hat{d}\theta_4
&~\in\Hom(\wedge^3\g_0\wedge\g_1,V_1),\\
(\D^\mathfrak{h}\theta_\varepsilon)_4&=&d_{l_3}\theta_0+d^{(1,0)}_{\phi}\theta_2+\widehat{d^\mathbb{V}}\theta_4&~\in\Hom(\wedge^4\g_0,V_0),\\
(\D^\mathfrak{h}\theta_\varepsilon)_5&=&d_{\phi_2}\theta_2+d_{l_3}\theta_3+d^{(1,0)}_{\phi}\theta_4&~\in\Hom(\wedge^5\g_0,V_1),
\end{array}\right.
  \end{equation*}
where $\theta_j=\theta_{\varepsilon
j}=(s^*\D^\g(\pi^*\lambda_\varepsilon))_j$ as in (\ref{3-ex cochain}).
By direct calculations, we have
\begin{eqnarray*}
(\D^\mathfrak{h}\theta_\varepsilon)_0(a,b)&=&-s_0(da)\triangleright{\lambda_0(b)}+\lambda_0[da,b]-\lambda_1(da,db)+d^\m
\lambda_2(da,b)\\
&&-s_0(db)\triangleright{\lambda_0(a)}+\lambda_0[db,a]-\lambda_1(db,da)+d^\m
\lambda_2(db,a)
\\ &&+d^{\mathbb{V}}\big(s_1(a)\triangleright \lambda_0(b)+s_1(b)\triangleright \lambda_0(a)-\lambda_2(da,b)+\lambda_2(a,db)\big)\\ &=&
\varphi_0\lambda_0(a)\triangleright
\lambda_0(b)+\varphi_0\lambda_0(b)\triangleright \lambda_0(a),
\end{eqnarray*}
which vanishes due to condition $(i)$ of Definition \ref{defi:Lie 2
cm}. Similarly, we can deduce that
$\D^\mathfrak{h}\theta_\varepsilon=0$.\qed\vspace{3mm}

As has been already demonstrated, to each strong crossed module
$\varepsilon=(\m,\g,\phi,\varphi)\in
\mathcal{C}(\mathfrak{h},\mathbb{V})$, one can obtain a $3$-cocycle
$\theta_\varepsilon\in C^3(\mathfrak{h},\mathbb{V})$. Then, define $\mu(\varepsilon)=[\theta_\varepsilon].$

Step $2$: We shall prove the canonical property of the map $\mu$. Namely, $\mu$ is independent of the choices
made of sections $s$ and $q$. Here, for future reference, we also
prove that if there is a strong map $\varepsilon\rightarrow{\varepsilon'}$
in $\mathcal{C}(\mathfrak{h},V)$, then $\theta_\varepsilon$ equals
$\theta_{\varepsilon'}$ in $\Ha^{3}(\mathfrak{h},\mathbb{V})$.
\begin{lem}\label{lem:section s}
$\theta_\varepsilon$ is independent of the choice of section $s$.
\end{lem}
\pf Suppose that $\bar{s}=(\bar{s}_0,\bar{s}_1)$ is another section
of $\pi$ and let $\bar{\theta}_\varepsilon$ be the 3-cocycle defined
using $\bar{s}$ instead of $s$. We need to prove that
$\bar{\theta}_\varepsilon$ coincides with $\theta_\varepsilon$ in
$\Ha^3(\mathfrak{h},\mathbb{V})$.

Since $\bar{s}$ and $s$ are both sections of $\pi$, there exist two
linear maps $t_{i}: \mathfrak{h}_i\longrightarrow \m_i$ with
$s_{i}-\bar{s}_i=\varphi_{i}\circ t_i.$ Construct four maps as
follows, for any $x,y,z\in{\mathfrak{h}_0},a\in{\mathfrak{h}_1}$,
\begin{equation}\label{sub}
\left\{\begin{array}{rcll}
B_0(a)&=&d^\m t_1(a)-t_0(da),\\
B_1(x,y)&=&\bar{s}_0(x)\triangleright{t_0(y)}-s_0(y)\triangleright{t_0(x)}-t_0[x,y],\\
B_2(x,a)&=&\bar{s}_0(x)\triangleright{t_1(a)}-s_1(a)\triangleright{t_0(x)}-t_1[x,a],\\
B_3(x,y,z)&=&-(s_0(x),s_0(y))\triangleright
t_0(z)-(\bar{s}_0(x),\bar{s}_0(y))\triangleright
t_0(z)+(s_0(x),\bar{s}_0(y))\triangleright t_0(z)\\
&&-(\bar{s}_0(y),s_0(z))\triangleright
t_0(x)-(\bar{s}_0(z),s_0(x))\triangleright t_0(y).\end{array}\right.
\end{equation}
Since $\varphi_0d=d^m\varphi_1$, it is obvious that
\begin{eqnarray*}
\varphi_0(B_0(a))&=&ds_1(a)-d\bar{s}_1(a)-s_0(da)+\bar{s}_0(da)\\
&=&\varphi_0(\lambda_0(a)-\bar{\lambda}_0(a)),
\end{eqnarray*}
which implies that
\[\lambda_0-\bar{\lambda}_0-B_0\in{\Hom(\mathfrak{h}_{1},V_0)}.\]
Similarly, relying on $\Pi=Id+\varphi$ is a strong homomorphism, we
can deduce that
$\lambda-\bar{\lambda}-B\in C^2(\mathfrak{h},\mathbb{V})$, where
$B=\sum_{i=0}^{3}B_i$. Furthermore, we claim that
$$\theta_{\varepsilon}-\bar{\theta}_{\varepsilon}=\D^\mathfrak{h}(\lambda-\bar{\lambda}-B).$$
By straightforward computations, we have
\begin{eqnarray}\label{1}
\nonumber(\theta_{\varepsilon}-\bar{\theta}_{\varepsilon})_0(x,a)
&=&s_{0}(x)\triangleright(\lambda_0-\bar{\lambda}_0)(a)-(\lambda_0-\bar{\lambda}_0)[x,a]
+\varphi_0t_0(x)\triangleright\bar{\lambda}_0(a)\\
&&+(\lambda_1-\bar{\lambda}_1)(x,da)-d^\m
(\lambda_2-\bar{\lambda}_2)(x,a),
\end{eqnarray}
and
\begin{eqnarray}\label{2}
\nonumber(\D^\mathfrak{h}(\lambda-\bar{\lambda}-B))_0(x,a)
&=&s_{0}(x)\triangleright(\lambda_0-\bar{\lambda}_0-B_0)(a)-(\lambda_0-\bar{\lambda}_0-B_0)[x,a]\\
&& +(\lambda_1-\bar{\lambda}_1-B_1)(x,da)-d^\m
(\lambda_2-\bar{\lambda}_2-B_2)(x,a).
\end{eqnarray}
Then substituting $B_i$ by the right hand sides of (\ref{sub}) and
taking into account condition $(i)$ of Definition \ref{defi:Lie 2
cm}, we get
\begin{eqnarray*}
&&(\ref{1})-(\ref{2})\\ &=&s_0(x)\triangleright (d^\m
t_1(a)-t_0(da))-d^\m
t_1[x,a]+t_0[x,da]+(\bar{s}_0(da)-d\bar{s}_1(a))\triangleright
t_0(x)\\ &&+
\bar{s}_0(x)\triangleright{t_0(da)}-s_0(da)\triangleright{t_0(x)}-t_0[x,da]
-d^\m\big(\bar{s}_0(x)\triangleright{t_1(a)}-s_1(a)\triangleright{t_0(x)}-t_1[x,a]\big)\\
&=&\varphi_0 t_0(x)\triangleright d^\m t_1(a)+d^\m(\varphi_1
t_1(a)\triangleright t_0(x)) -\varphi_0 t_0(x)\triangleright
t_0(da)-\varphi_0 t_0(da)\triangleright t_0(x)\\ &=&0.
\end{eqnarray*}
Similarly, we can verify that
$\theta_{\varepsilon}-\bar{\theta}_{\varepsilon}=\D^\mathfrak{h}(\lambda-\bar{\lambda}-B)$,
which proves that the class of $\theta_\varepsilon$ is independent
of the section $s$. \qed\vspace{3mm}

\begin{lem}\label{map}
$\theta_\varepsilon$ is independent of the choice of section $q$ and
a strong map $\varepsilon\longrightarrow{\varepsilon'}.$
\end{lem}
\pf Consider a strong map $(F,G):\varepsilon\longrightarrow{\varepsilon'}$
as in Definition \ref{defi:el}. Let
$s:\mathfrak{h}\longrightarrow{\mathfrak{g}}$ and
$q:\Img\varphi\longrightarrow{\mathfrak{m}}$ be sections of $\pi$
and $\varphi$ and let
$s':\mathfrak{h'}\longrightarrow{\mathfrak{g'}}$ and
$q':\Img\varphi'\longrightarrow{\mathfrak{m'}}$ be sections of
$\pi'$ and $\varphi'$. Since $\pi'(Gs)=\pi(s)=Id$, we get another
section $Gs=(G_{1}s_1,G_{0}s_0)$ of $\pi'$. Taking into account
Lemma \ref{lem:section s}, we choose $s'=Gs$. Set
\begin{equation*}
\left\{\begin{array}{rcl}
B_0(a)&=&(F_{0}q_0-q'_0G_{0})(ds_1(a)-s_{0}d(a)),\\
B_1(x,y)&=&(F_{0}q_{0}-q'_{0}G_0)([s_0(x),s_0(y)]-s_0[x,y]),\\
B_2(x,a)&=&(F_{1}q_1-q'_{1}G_{1})([s_0(x),s_1(a)]-s_1[x,a]),\\
B_3(x,y,z)&=&(F_{1}q_1-q'_{1}G_{1})(l_3(s_0(x),s_0(y),s_0(z))-s_{1}l_3(x,y,z)).\end{array}\right.
\end{equation*}
Noticing that $\varphi'F=G\varphi$, it is obvious that
$B=\sum_{i=0}^3 B_i\in C^2(\mathfrak{h},\mathbb{V}).$ Furthermore,
relying on the properties of $(F,G)$, we obtain
\begin{eqnarray*}
&&(\theta_{\varepsilon0}-\theta_{\varepsilon'0})(x,a)\\
&=&F_0(\theta_{\varepsilon0}(x,a))-\theta_{\varepsilon'0}(x,a)\\ &=&
F_0\big(s_0(x)\triangleright{q_0(ds_1(a)-s_0(da))}-
q_0(ds_1[x,a]-s_0d[x,a])\\ &&+q_0([s_0(x),s_0(da)]-s_0[x,da])-d^\m
q_1([s_0(x),s_1(a)]-s_1[x,a])\big)\\
&&-G_0s_0(x)\triangleright'{q'_0(d'G_1s_1(a)-G_0s_0(da))}+
q'_0(d'G_1s_1[x,a]-G_0s_0d[x,a])\\
&&-q'_0([G_0s_0(x),G_0s_0(da)]-G_0s_0[x,da])+d^{\m'}q'_1([G_0s_0(x),G_1s_1(a)]-G_1s_1[x,a])
\\ &=&G_{0}s_0(x)\triangleright'(F_{0}q_0-q'_0G_{0})(ds_1(a)-s_0(da))
-(F_{0}q_0-q'_0G_{0})(ds_1[x,a]-s_{0}d[x,a])\\
&&+(F_{0}q_{0}-q'_{0}G_0)([s_0(x),s_0(da)]-s_0[x,da])-d^{\m'}
(F_{1}q_1-q'_{1}G_{1})([s_0(x),s_1(a)]-s_1[x,a]))\\
&=&(d^{(1,0)}_{\phi'}B_0+\hat{d}B_1+\widehat{d^{\m'}}B_2)(x,a).
\end{eqnarray*}
Similarly, considering the fact that $\phi$ and $\phi'$ induce the
same $\mathfrak{h}$-module structure on $\mathbb{V}$, we get
$$\theta_\varepsilon-\theta_{\varepsilon'}=\D^{\mathfrak{h}}B.$$
 This finishes the proof. \qed\vspace{3mm}

Step $3$: We show $\mu$ is a surjection, which follows from Corollary \ref{main2} and the following lemma.
\begin{lem}\label{sur}
For any $[\theta]\in \Ha^3(\mathfrak{h},\mathbb{V})$, we have $\mu(\mathcal{C}(\mathfrak{h},\mathbb{V})_{[\theta]})=[\theta]$.
\end{lem}
 \pf For a
crossed module $\varepsilon_\lambda=(\ker
\pi\oplus_\lambda\mathbb{V},\mathfrak{F}(\mathfrak{h}),\phi^\lambda,\varphi)
\in\mathcal{C}(\mathfrak{h},\mathbb{V})_{[\theta]},$
consider the complex $$\varepsilon_\lambda:
0\rightarrow\mathbb{V}\stackrel{i}\longrightarrow\ker
\pi\oplus_\lambda\mathbb{V}\stackrel{\varphi}\longrightarrow\mathfrak{F}(\mathfrak{h})\stackrel{\pi}\longrightarrow\mathfrak{h}\rightarrow0.$$
Choosing any section $s$ of $\pi$ and defining section $q$ of $\varphi$ on
$\Img \varphi$ by $q(\alpha)=(\alpha,0), \forall\alpha\in
\ker\pi$, we get a $3$-cocycle $\theta_{\varepsilon_\lambda}\in
C^3(\mathfrak{h},\mathbb{V})$. We
claim that
$\theta_{\varepsilon_\lambda}+\D^{\mathfrak{h}}s^*\lambda=\theta$,
which implies that
$\mu(\varepsilon_\lambda)=[\theta_{\varepsilon_\lambda}]=[\theta]$.
Actually, by (\ref{34}) and (\ref{3-ex cochain}), we have
\begin{eqnarray*}
(\theta_{\varepsilon_\lambda})_0(x,a)&=&s_0x\triangleright_\lambda(ds_1a-s_0da)
-ds_1[x,a] +s_0d[x,a]\\ &&+[s_0x,s_0da]
-s_0[x,da]-d^\lambda([s_0x,s_1a]-s_1[x,a])
\\ &=&\lambda_1(s_0x,ds_1a-s_0da)
-\lambda_0([s_0x,s_1a]-s_1[x,a]),
\end{eqnarray*}
and
\begin{eqnarray*}
(\D^{\mathfrak{h}}s^*\lambda)_0(x,a)
&=&
x\triangleright
\lambda_0(s_1a)-\lambda_0s_1[x,a]+\lambda_1(s_0x,s_0da)-d^\mathbb{V}\lambda_2(s_0x,s_1a).
\end{eqnarray*}
Adding them together, we get
\begin{eqnarray*}
(\theta_{\varepsilon_\lambda}+\D^{\mathfrak{h}}s^*\lambda)_0(x,a)&=&x\triangleright
\lambda_0(s_1a)-\lambda_0[s_0x,s_1a]+\lambda_1(s_0x,ds_1a)
-d^\mathbb{V}\lambda_2(s_0x,s_1a)\\
&=&(\D^{\mathfrak{F}(\mathfrak{h})}\lambda)_0(s_0x,s_1a)=
(\pi^*\theta)_0(s_0x,s_1a)\\ &=&\theta_0(x,a).
\end{eqnarray*}
Likewise, we can obtain
$\theta_{\varepsilon_\lambda}+\D^{\mathfrak{h}}s^*\lambda=\theta$.
This finishes the proof. \qed\vspace{3mm}

Step $4$: We shall prove for two crossed modules $\varepsilon,\varepsilon'\in\mathcal{C}(\mathfrak{h},\mathbb{V})$,
$\mu(\varepsilon)=\mu(\varepsilon')$ iff $\varepsilon\sim\varepsilon'$. Then, the map $\mu: \mathcal{C}(\mathfrak{h},\mathbb{V})\longrightarrow\Ha^{3}(\mathfrak{h},\mathbb{V})$ induces a bijection between $\mathcal{C}(\mathfrak{h},\mathbb{V})/_\sim$ and $\Ha^{3}(\mathfrak{h},\mathbb{V})$.

A direct consequence of Lemma \ref{map} and \ref{sur} is:
\begin{cor}\label{can}
If $\varepsilon\sim\varepsilon'$, we have $\mu(\varepsilon)=\mu(\varepsilon')$.
\end{cor}
\begin{pro}\label{inj}For a $3$-cocycle $\theta\in C^3(\mathfrak{h},\mathbb{V})$,
$\mu(\varepsilon)=[\theta]$ if and only if there exists a crossed
module $\varepsilon_\lambda\in
\mathcal{C}(\mathfrak{h},\mathbb{V})_{[\theta]}$ and a strong map
$(F,G):\varepsilon_\lambda\longrightarrow \varepsilon$. That is, if $\mu(\varepsilon)=\mu(\varepsilon')$, then $\varepsilon\sim\varepsilon'$.
\end{pro}
\pf For such a crossed module
$\varepsilon=(\m,\g,\phi,\varphi)\in\mathcal{C}(\mathfrak{h},\mathbb{V})$,
choosing sections $s:\mathfrak{h}\longrightarrow \g$ and
$q:\Img\varphi\longrightarrow \m$ of $\pi'$ and $\varphi$
respectively, we can construct $\lambda_\varepsilon$ and then a
$3$-cocycle $\theta_\varepsilon \in C^3(\mathfrak{h},\mathbb{V})$ as in (\ref{34}) such that $[\theta_\varepsilon]=[\theta]$.
Let $G:\mathfrak{F}(\mathfrak{h})\longrightarrow \g$ be the strong
Lie $2$-algebra homomorphism induced by $s$ (the property of free
Lie $2$-algebras). Define linear maps
$\psi=(\psi_0,\psi_1):\mathfrak{F}(\mathfrak{h})\longrightarrow \m$
by $\psi(\bar{x})=q(G\bar{x}-s\pi\bar{x}), \forall \bar{x}\in
\mathfrak{F}(\mathfrak{h}).$ Next, construct four maps:
\begin{equation*}
\left\{\begin{array}{rcl}
\lambda_0(\bar{a})&=&\lambda_{\varepsilon0}(\pi_1\bar{a})+d^{\m}\psi_1 \bar{a}-\psi_0 d\bar{a},\\
\lambda_1(\bar{x},\bar{y})&=&\lambda_{\varepsilon1}(\pi_0\bar{x},\pi_0\bar{y})
-\psi_0[\bar{x},\bar{y}]-[\psi_0\bar{x},\psi_0\bar{y}]+G_0\bar{x}\triangleright\psi_0\bar{y}
-G_0\bar{y}\triangleright\psi_0\bar{x},\\
\lambda_2(\bar{x},\bar{a})&=&\lambda_{\varepsilon2}(\pi_1\bar{x},\pi_1\bar{a})
-\psi_1[\bar{x},\bar{a}]-[\psi_0\bar{x},\psi_1\bar{a}]+G_0\bar{x}\triangleright\psi_1\bar{a}
-G_1\bar{a}\triangleright\psi_0\bar{x},\\
\lambda_3(\bar{x},\bar{y},\bar{z})&=&\lambda_{\varepsilon3}(\pi_0\bar{x},\pi_0\bar{y},\pi_0\bar{z})-\psi_1l_3(\bar{x},\bar{y},\bar{z})+
l^\m_3(\psi_0\bar{x},\psi_0\bar{y},\psi_0\bar{z})\\
&&-(l_{\phi_0(G\bar{x})}(\psi_0\bar{y},\psi_0\bar{z})+(G_0\bar{x},
G_0\bar{y})\triangleright \psi_0\bar{z}+c.p.),\end{array}\right.
\end{equation*}
for any $\bar{x},\bar{y},\bar{z}\in \mathfrak{F}(\mathfrak{h})_0, \bar{a}\in \mathfrak{F}(\mathfrak{h})_1.$ Since 
$\pi, G$ and $\varphi$ commute with $d$, we have
\begin{eqnarray*}
\varphi_0(\lambda_0(\bar{a}))&=&ds_1\pi_1\bar{a}-s_0d\pi_1\bar{a}+d(G_1\bar{a}-s_1\pi_1\bar{a})-(G_0d\bar{a}-
s_0\pi_0d\bar{a})\\ &=&0.
\end{eqnarray*}
Similarly, we can get $\varphi\circ \lambda_i=0, i=1,2,3$, that is,
 $\lambda=\sum_{i=0}^{3} \lambda_i\in C^2(\mathfrak{F}(\mathfrak{h}),\mathbb{V}).$

Moreover, we claim that
$\D^{\mathfrak{F}(\mathfrak{h})}\lambda=\pi^*\theta_\varepsilon$. By
straightforward calculations, we get
\begin{eqnarray*}
(\D^{\mathfrak{F}(\mathfrak{h})}\lambda)_0(\bar{x},\bar{a})&=&\pi_0\bar{x}\triangleright
\lambda_0(\bar{a})-\lambda_0[\bar{x},\bar{a}]+\lambda_1(\bar{x},d\bar{a})
-d^\m\lambda_2(\bar{x},\bar{a})\\
&=&\pi_0\bar{x}\triangleright\big(\lambda_{\varepsilon0}(\pi_1\bar{a})+d^{\m}\psi_1
\bar{a}-\psi_0
d\bar{a}\big)-\lambda_{\varepsilon0}(\pi_1[\bar{x},\bar{a}])-d^{\m}\psi_1
[\bar{x},\bar{a}]+\psi_0 d[\bar{x},\bar{a}]\\ &&
+\lambda_{\varepsilon1}(\pi_0\bar{x},\pi_0 d\bar{a})
-\psi_0[\bar{x},d\bar{a}]-[\psi_0\bar{x},\psi_0d\bar{a}]+G_0\bar{x}\triangleright\psi_0d\bar{a}
-G_0d\bar{a}\triangleright\psi_0\bar{x}\\
&&-d^\m\big(\lambda_{\varepsilon2}(\pi_0\bar{x},\pi_1\bar{a})
-\psi_1[\bar{x},\bar{a}]-[\psi_0\bar{x},\psi_1\bar{a}]+G_0\bar{x}\triangleright\psi_1\bar{a}
-G_1\bar{a}\triangleright\psi_0\bar{x}\big)\\
&=&s_0\pi_0\bar{x}\triangleright\lambda_{\varepsilon0}(\pi_1\bar{a})
-\lambda_{\varepsilon0}(\pi_1[\bar{x},\bar{a}])+\lambda_{\varepsilon1}(\pi_0\bar{x},\pi_0
d\bar{a})-d^\m\lambda_{\varepsilon2}(\pi_0\bar{x},\pi_1\bar{a})\\
&=&\theta_{\varepsilon0}(\pi_0\bar{x},\pi_1\bar{a})\\
&=&(\pi^*\theta_{\varepsilon0})(\bar{x},\bar{a}),
\end{eqnarray*}
where we have used  condition $(i)$ of Definition \ref{defi:Lie 2
cm}. Likewise, we deduce that
$\D^{\mathfrak{F}(\mathfrak{h})}\lambda=\pi^*\theta_\varepsilon$.

So we can use such a defined $\lambda$ to construct a crossed module
$\varepsilon_{\lambda}$ as in Proposition \ref{pro:surjection}. In the following, we prove
 there is a map $(F,G)$ from $\varepsilon_{\lambda}$ to
$\varepsilon$ in $\mathcal{C}(\mathfrak{h},\mathbb{V})$:
\begin{equation*}
\CD
  0 @> >>  \mathbb{V} @>i>> \ker \pi\oplus_{\lambda}\mathbb{V} @>\varphi^\lambda>> \mathfrak{F}(\mathfrak{h})@>\pi>> \mathfrak{h} @> >> 0 \\
  @. @V Id VV @V FVV @VG VV @V Id VV @.  \\
  0 @> >> \mathbb{V} @>i'>> \mathfrak{m} @>\varphi>> \mathfrak{g}@>\pi'>> \mathfrak{h} @> >> 0,
\endCD
\end{equation*}
where $F(\alpha+u)=i'(u)+qG(\alpha), \forall \alpha\in
\ker\pi,u\in\mathbb{V}.$ Indeed, we shall prove $F$ is a strong Lie
$2$-algebra homomorphism compatible with the actions, that is, for
any $\alpha\in \ker \pi,u\in\mathbb{V}, \bar{x}\in
\mathfrak{F}(\mathfrak{h}),$
 $$F(\bar{x}\triangleright_{\lambda}(\alpha+u))=G(\bar{x})\triangleright F(\alpha+u), \ \ \ F((\bar{x},\bar{y})\triangleright_{\lambda}(\alpha+u))=(G\bar{x},G\bar{y})\triangleright F(\alpha+u). $$
 We only sketch the proof of compatibility. For any $\alpha\in \ker\pi_0,u\in V_0$ and $\bar{x}_0\in
\mathfrak{F}(\mathfrak{h})_0$, we have
\begin{eqnarray*}
F_0(\bar{x}\triangleright_{\lambda}(\alpha+u))&=&F_0([\bar{x},\alpha]+\lambda_1(\bar{x},\alpha)+\pi_0\bar{x}\triangleright
u)\\ &=&\lambda_1(\bar{x},\alpha)+\pi_0\bar{x}\triangleright
i'_0u+q_0G_0[\bar{x},\alpha]\\ &=&
-q_0G_0[\bar{x},\alpha]-[q_0G_0\bar{x}-q_0s_0\pi_0\bar{x},q_0G_0\alpha]+G_0\bar{x}\triangleright
q_0G_0 \alpha\\ &&-G_0\alpha\triangleright (q_0G_0
\bar{x}-q_0s_0\pi_0\bar{x})+\pi_0\bar{x}\triangleright
i'_0u+q_0G_0[\bar{x},\alpha]\\ &=&\pi_0\bar{x}\triangleright i'_0u+
G_0\bar{x}\triangleright q_0G_0 \alpha\\ &=&G_0\bar{x}\triangleright
F_0(\alpha+u),
\end{eqnarray*}
where we have used the condition $(i)$ of Definition \ref{defi:Lie 2
cm} and $\pi_0\bar{x}\triangleright
i'_0u=\pi'_0G_0\bar{x}\triangleright i'_0u=G_0\bar{x}\triangleright
i'_0u.$ Next, since $(\pi_0\bar{x},\pi_0\bar{y})\triangleright
i'_0u=(G_0\bar{x},G_0\bar{y})\triangleright i'_0u$ due to the fixed
action of $\mathfrak{h}$ on $\mathbb{V}$, we have
\begin{eqnarray*}
&&F_1((\bar{x},\bar{y})\triangleright_{\lambda}(\alpha+u))\\
&=&F_1(l_3(\bar{x},\bar{y},\alpha)+\lambda_3(\bar{x},\bar{y},\alpha)+(\pi_0\bar{x},
\pi_0\bar{y})\triangleright u)\\ &=&
-q_1G_1l_3(\bar{x},\bar{y},\alpha)+l^{\m}_3(q_0G_0\bar{x}-q_0s_0\pi_0\bar{x},q_0G_0\bar{y}-q_0s_0\pi_0\bar{y},
q_0G_0\alpha)\\ &&-(G_0\bar{x},G_0\bar{y})\triangleright
q_0G_0\alpha-(G_0\bar{y},G_0\alpha)\triangleright(q_0G_0\bar{x}-q_0s_0\pi_0\bar{x})
-(G_0z,G_0\bar{x})\triangleright(q_0G_0\bar{y}-q_0s_0\pi_0\bar{y})\\
&&+(G_0\bar{x},G_0\bar{y}-s_0\pi_0\bar{y})\triangleright
q_0G_0\alpha+(G_0\bar{y},G_0\alpha)\triangleright
(q_0G_0\bar{x}-q_0s_0\pi_0\bar{x})\\
&&+(G_0\alpha,G_0\bar{x}-s_0\pi_0\bar{x})\triangleright
(q_0G_0\bar{y}-q_0s_0\pi_0\bar{y})+q_1G_1l_3(\bar{x},\bar{y},\alpha)+(\pi_0\bar{x},
\pi_0\bar{y})\triangleright i'_0u\\
&=&(\pi_0\bar{x},\pi_0\bar{y})\triangleright
i'_0u-(G_0\bar{x},G_0\bar{y})\triangleright q_0G_0\alpha\\
&=&(G_0\bar{x},G_0\bar{y})\triangleright F_0(\alpha+u),
\end{eqnarray*}
where we have used condition $(i),(ii),(iii)$ of Definition
\ref{defi:Lie 2 cm}. The other hand is a consequence of Lemma
\ref{map} and \ref{sur}. This ends the
proof.\qed\vspace{3mm}

Next, we give an alternative description of the equivalence relation in Theorem \ref{thm:main}, which is similar to the statements in \cite{Wagemann, Casas}.
Two strong crossed modules $(\m,\g,\phi,\varphi),
(\m',\g',\phi',\varphi')\in \mathcal{C}(\mathfrak{h},\mathbb{V})$
are called {\bf elementary equivalent} if there is a morphism of
crossed modules $(F,G,\tau)$ such that
the diagram
\begin{equation*}
\CD
  0 @> >>  \mathbb{V} @>i>> \mathfrak{m} @>\varphi>> \mathfrak{g}@>\pi>> \mathfrak{h} @> >> 0 \\
 @. @V Id VV @V FVV @VG VV @V Id VV @.  \\
 0 @> >> \mathbb{V} @>i'>> \mathfrak{m'} @>\varphi'>> \mathfrak{g'}@>\pi'>> \mathfrak{h} @> >> 0,
\endCD
\end{equation*}
is commutative and
$$G_2=0; \, \ \ \Img\tau\subset i'(\mathbb{V});\ \ \
\tau(\g_0\wedge i(\mathbb{V}))=0;\ \ \
\tau(\varphi_0\alpha,\beta)=\tau(\alpha,\varphi_0\beta),\ \ \
\forall \alpha,\beta\in \m_0.$$
By a straightforward but tedious computation combined with the proof of Theorem \ref{thm:main}, we obtain the following proposition.
\begin{pro}
The equivalence relation generated by elementary equivalence relation coincides with the
equivalence relation $\sim$ in Theorem \ref{thm:main}.
\end{pro}
One key step in the proof of Theorem \ref{thm:main} is the construction of the map $\mu: \mathcal{C}(\mathfrak{h},\mathbb{V})\rightarrow \Ha^3(\mathfrak{h},\mathbb{V})$. The following example illustrates that for a particular class of crossed modules given in Example \ref{2-vector spaces}, the map $\mu$ is closely related to the connecting map in the long exact sequence of cohomology groups.
\begin{ex}{\rm Consider the strong crossed module $\varepsilon=(\mathbb{I},\mathfrak{h}\oplus_\lambda \mathbb{Q},\phi,\varphi)\in \mathcal{C}(\mathfrak{h},\mathbb{V})$ obtained in Example \ref{2-vector spaces}.
The sequence (\ref{ex se1}) induces a short sequence of complexes
$$ 0\rightarrow C^*(\mathfrak{h},\mathbb{V})\stackrel{\bar{p}}\longrightarrow C^*(\mathfrak{h},\mathbb{I})\stackrel{\bar{q}}\longrightarrow C^*(\mathfrak{h},\mathbb{Q})\rightarrow0.$$
Since for arbitrary $\mathfrak{h}$-module homomorphism $f$, we have
$\bar{f}\circ \D^\mathfrak{h}=\D^\mathfrak{h}\circ \bar{f}$. Thus
$p,q$ induce maps between cohomology groups. Moreover, similar to
the process in homological algebra, we can construct a connecting
homomorphism $\partial: \Ha^*(\mathfrak{h},\mathbb{Q})\rightarrow
\Ha^{*+1}(\mathfrak{h},\mathbb{V})$ such that
$$ \cdots \rightarrow \Ha^*(\mathfrak{h},\mathbb{V})\stackrel{\bar{p}}\longrightarrow \Ha^*(\mathfrak{h},\mathbb{I})\stackrel{\bar{q}}\longrightarrow \Ha^*(\mathfrak{h},\mathbb{Q})\stackrel{\partial}\longrightarrow \Ha^{*+1}(\mathfrak{h},\mathbb{V})\rightarrow \cdots$$ is a long exact sequence of cohomology groups.

In particular, we have $\mu(\varepsilon)=\partial[\lambda]$. See
\cite[Theorem 3]{Wagemann} for more details.}
\end{ex}

 \end{document}